\documentclass[11pt]{amsart}

\usepackage[leqno]{amsmath}
\usepackage{amssymb}
\usepackage{amsfonts}
\usepackage{amsxtra}
\usepackage{graphicx, color, epsfig}

\swapnumbers

\newcommand{\pairing}[2]{\langle\,#1\,,\,#2\,\rangle}

\newcommand{\lk}[1]{\mathcal{#1}}

\newcommand{\F}{\mathbb{F}}
\newcommand{\Z}{\mathbb{Z}}

\newcommand{\R}{\mathbb{R}}

\newcommand{\area}[1]{[#1]}
\newcommand{\ring}[1]{\mathbb{P}_{#1}}
\newcommand{\field}[1]{\mathbb{F}_{#1}}
\newcommand{\cross}[1]{\ensuremath{\mathrm{\textsc{cr}}(#1)}}
\newcommand{\arcs}[1]{\ensuremath{\mathrm{\textsc{arc}}(#1)}}
\newcommand{\regions}[1]{\mathfrak{F}_{L}}
\newcommand{\graph}[1]{\Gamma_{L}}
\newcommand{\circles}[2][]{\mathrm{\textsc{cir}}_{#1}(#2)}
\newcommand{\glyph}[2]{\mathrm{\textsc{glyph}(#1,#2)}}
\newcommand{\gly}[1]{\mathrm{\textsc{gl}(#1)}}
\newcommand{\res}[1]{\ensuremath{\mathrm{\textsc{res}}(#1)}}

\newcommand{\lra}{\longrightarrow}

\newcommand{\inlinediag}[2][0.33]{\includegraphics[scale=#1]{#2}}
\newcommand{\diagram}[4][1.0]{\begin{center}\begin{figure}\includegraphics[scale=#1]{#3}\caption{#4}\label{fig:#2}\end{figure}\end{center}}

\newtheorem{theorem}{Theorem}[section]
\newtheorem{defn}[theorem]{Definition}
\newtheorem{lemma}[theorem]{Lemma}
\newtheorem{cor}[theorem]{Corollary}
\newtheorem{prop}[theorem]{Proposition}

\newcommand{\POz}{P.\ Ozsv{\'a}th\,}
\newcommand{\ZSz}{Z.\ Szab{\'o}\,}

\newcommand{\APS}{M. Asaeda, J. Przytycki, and A. Sikora\,}

\title{Twisted Skein Homology}
\author{Nguyen D. Duong \& Lawrence P. Roberts}

\begin{document}
\maketitle


\section{Introduction}

\noindent In \cite{LR} the second author described how to take the diagram $L$ of a link $\lk{L} \subset S^{3}$ and modify  the reduced, characteristic 2, $\delta$-graded Khovanov chain complex for $L$ to get an explicit knot homology theory built out of the spanning trees of a Tait graph for $L$. This occurs by working through an intermediary called ``totally twisted Khovanov homology''. In \cite{TJ} Thomas Jaeger proved that for a knot $K$ the spanning tree chain complex computes the reduced Khovanov homology of $K$ (over a specific field). In addition, Andrew Manion has extended the totally twisted theory to odd Khovanov homology, and thus characteristic 0 coefficients, \cite{Mann}. Furthermore, other knot homology theories also admit an explicit description through spanning trees. J. Baldwin and A. Levine describe such a complex for knot Floer homology in \cite{Bald}, while another theory has been discovered by \POz, \ZSz, explored by J. Baldwin, and described by Daniel and Igor Kriz in \cite{Kriz}.\\
\ \\
\noindent In an effort to expand the reach of such theories, the authors took the known extensions of Khovanov homology given by
\APS in \cite{APS1}, \cite{APS2} and twisted them as well. The resulting boundary map is more complicated than in \cite{LR}, but still gives a theory that is invariant under link isotopies. \\
\ \\
\noindent To describe the chain complex, first recall the setting of \cite{APS1}, \cite{APS2}. Let $\Sigma$ be a compact, {\em oriented} surface and let $B \subset \partial \Sigma$ be a fixed collection of $2l$-points. An unoriented {\em tangle} $\mathcal{T}$ for $(\Sigma, B)$ is a finite collection of $l$ arcs $\alpha_{i}$ and some number of circles $\gamma_{i}$ which are disjointly embedded in $\Sigma \times (-1,1)$, such that the boundaries of the arcs exhaust $B \times \{0\}$. When $B = \emptyset$, we will also call $\mathcal{T}$ a link in $\Sigma \times (-1,1)$. Two such tangles will be equivalent if they are related by boundary preserving ambient isotopy. \\
\ \\
\noindent $\mathcal{T}$ can be isotoped to have a generic image under the projection $\pi: \Sigma \times (-1,1) \rightarrow \Sigma$, allowing us to study $\mathcal{T}$ through its crossing diagrams in $\Sigma$. We use such a diagram $T$ to define
a chain complex generated by the resolutions of $T$ which consist solely of noncontractible circles and arcs in $\Sigma$. \\
\ \\
\noindent More specifically, we label each edge of $T$ with a distinct formal variable $x_{f}$. Let $\field{T}$ be
the field of rational functions $\Z/2\Z(x_{1}, \ldots, x_{n})$. In this paper, we define chain complexes whose generators consist of resolutions of the diagram $T$:

\noindent 

\begin{defn}
For each subset $S$ of the crossings in $T$ the {\em resolution} $r_{S}$ of $T$ is the surface tangle in $\Sigma$, considered up to isotopy, found by locally replacing each crossing $c \in \cross{T}$ according to the rule
$$
\inlinediag{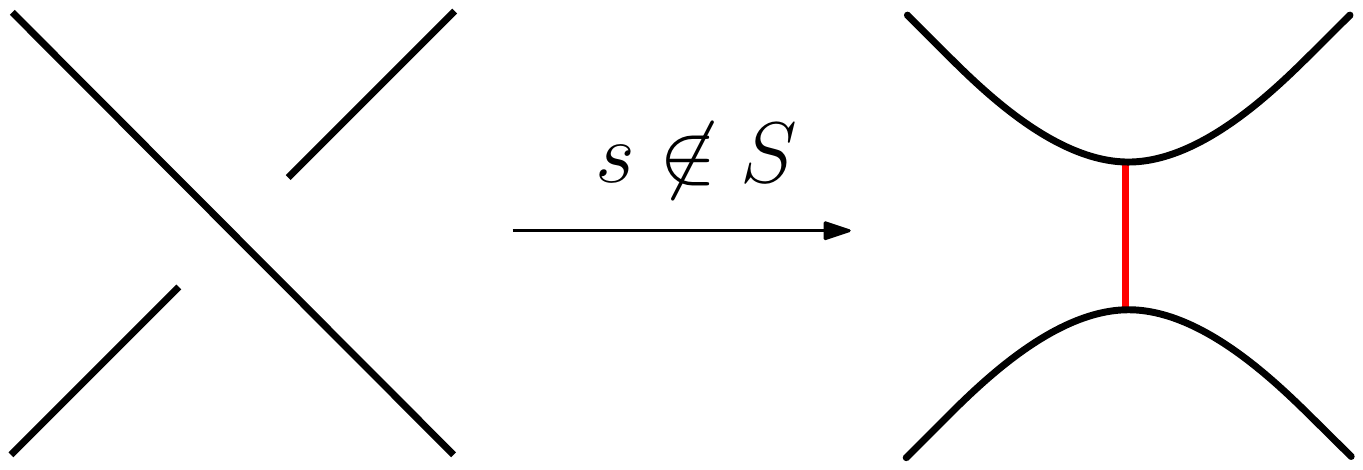} \hspace{.5in} \inlinediag{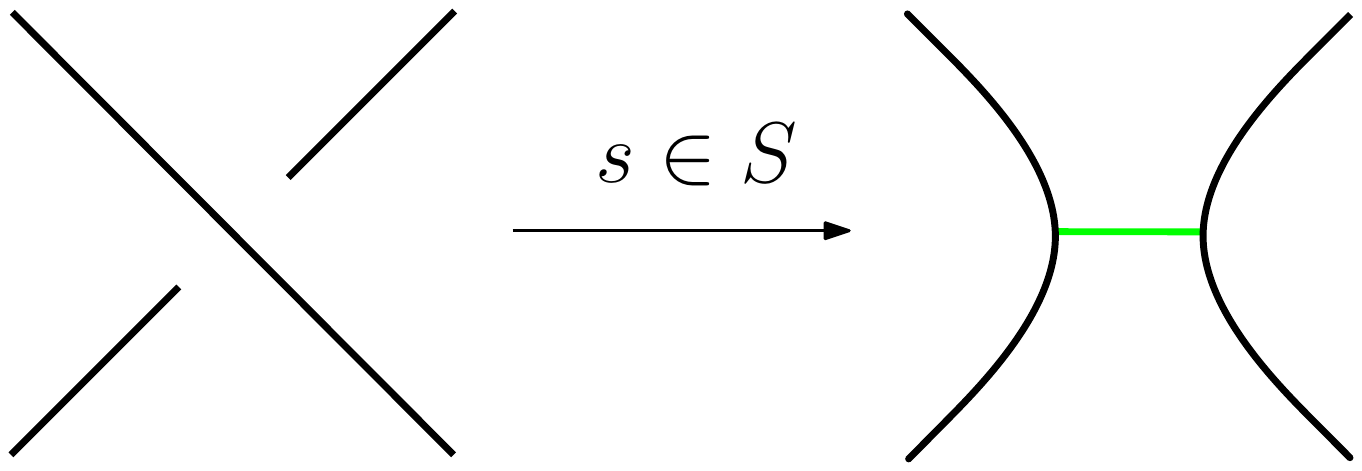} 
$$ 
The set of resolutions for $T$ will be denoted $\res{T}$. For each resolution $r$, denote by $|r|$ the number of elements in the corresponding subset $S \subset \cross{T}$. 
\end{defn}

\noindent Each resolution $r$ consists of a collection of $l$-arcs, some contractible circles, and come noncontractible circles, each with segments labeled by the $x_{f}$. Let $\lk{O}_{n}(T)$ be the set of pairs $(r, s)$  of a resolution $r$ consisting solely of arcs and noncontractible circles and a map $s$ assigning a number $\{+1,-1\}$ to each noncontractible circle in $r$. We will also call such a pair $r$  when it is clear from the context that an assignment is given. To such an $r$ 
we associate a glyph:  

\begin{defn}
A {\em glyph} for  $(\Sigma, B)$ is an equivalence class of disjoint proper embeddings in $\Sigma$ of
\begin{enumerate}
\item a finite collection of noncontractible circles $C_{i}$, each labeled with an integer $n_{i} \in \Z$, and
\item a finite collection of arcs, whose boundaries exhaust $B$   
\end{enumerate}
considered up to the equivalence generated by boundary preserving surface isotopy, and the erasure/introduction of $0$ labeled noncontractible circles. We will denote the set of glyphs by $\glyph{\Sigma}{B}$. In $\glyph{\Sigma}{\emptyset}$ there is a {\em trivial glyph}, the class where the number of circles embedded is $0$. 
\end{defn} 

\noindent The glyph $\gly{r}$ for $r$ is found by identifying the noncontractible circles in $r$ which cobound annuli, and adding their decorations, to get a single circle labeled with an integer. Given a glyph $b \in \glyph{\Sigma}{B}$ let
$$
\lk{O}_{n,i}(T,b) = \big\{\,r \in \lk{O}_{n}(T)\, \big|\,|r|=i,\,\gly{r}=b\,\big\}
$$
We use the elements of $\lk{O}_{n,i}(T,b)$ to define the $i^{th}$ grading in a chain complex:
$$
\widetilde{CT}_{i}(T,b):= Span_{\field{T}} \big \{ r \in \lk O_{n,i} (T,b)  \big \} 
$$
where the boundary map $\partial_{n,i}: \widetilde{CT}_{i}(T,b) \rightarrow \widetilde{CT}_{i+2}(T,b)$ is given by
$$
\partial_{n,i}(r)= \sum_{r' \in  \lk{O}_{n,i+2}(T,b)} \pairing{r}{r'}\,r'
$$ 
The coefficient $\pairing{r}{r'}$ is calculated from the weights assigned to the edges. First, $\pairing{r}{r'} = 0$ unless there are two crossings $c_{1}$ and $c_{2}$ which are resolved with the $0$ resolution in $r$ and the $1$ resolution in $r'$. Let $r_{01}$ be the resolution where $c_{1}$ is resolved with the $0$ resolution, and $c_{2}$ resolved with the $1$ resolution. Furthermore, if $C$ is a {\em contractible} circle,
let $[C] = \sum_{f \in \arcs{C}} x_{f}$. Then 

\begin{defn}
Suppose $r \in  \lk{O}_{n,i}(T,b)$ and $r' \in \lk{O}_{n,i+2}(T,b)$, and the resolutions $r_{01}, r_{10}$ defined above have $\gly{r_{01}}=\gly{r_{10}}=b$. Then
\begin{equation}
\pairing{r}{r'} = \left\{\begin{array}{cl} 
\displaystyle{\frac{1}{[C_{10}]} + \frac{1}{[C_{01}]}} & r_{10}\mathrm{\ and\ }r_{01}\mathrm{\ have\ contractible\ circles\ }C_{10}\mathrm{\ and\ }C_{01}\mathrm{\ respectively} \\
\displaystyle{\frac{1}{[C_{10}]}} & r_{10}\mathrm{\ has\ a\ contractible\ circle\ }C_{10}\mathrm{\ but\ }r_{01}\mathrm{\ does\ not}\\
\displaystyle{\frac{1}{[C_{01}]}} & r_{01}\mathrm{\ has\ a\ contractible\ circle\ }C_{01}\mathrm{\ but\ }r_{10}\mathrm{\ does\ not}\\
0 & \mathrm{\ neither\ }r_{10}\mathrm{\ nor\ }r_{01}\mathrm{\ has\ a\ contractible\ circle}\\
\end{array}\right.
\end{equation}
If $r'$ is a resolution with $|r'| = |r| + 2$ but no pair of crossings $c_{1}$ and $c_{2}$ exists with $r' = r_{11}$ as above then $\pairing{r}{r'} = 0$.
\end{defn}

\noindent For each glyph $b$, $\partial_{n,\ast}$ is a boundary map, and $\widetilde{CT}_{\ast}(T,b)$ is a chain complex. If $\mathcal{T}$ is oriented, let $n_{+}$ be the number of positive crossings in $T$ (relative to the orientation on $\Sigma$). Then

\begin{theorem}
For each glyph $b \in \glyph{\Sigma}{B}$, the (stable) chain homotopy type of $(\widetilde{CT}_{*}(T,b), \partial_{n})[-n_{+}(T)]$ is invariant under the Reidemeister I,II, and III moves in $\Sigma$, and thus is an invariant of $\lk{T}$. 
\end{theorem}

\noindent {\bf Note:} Stable equivalence is a technical means for relating the complexes despite that $\field{T}$ can change
under the Reidemeister moves. See \cite{LR} for more details. The notation $[-n_{+}(T)]$ indicates a shift in the gradings which is described more fully in the next section.\\

\begin{defn}
The homology groups of $(CT_{*}(T,b), \partial_{n})[-n_{+}(T)]$ will be denoted $\underline{HT}_{*}(\mathcal{T},b)$ and 
$$\underline{HT}_{*}(\mathcal{T})  = \displaystyle{\bigoplus_{b \in \glyph{\Sigma}{B}} \underline{HT}_{*}(\mathcal{T},b)}$$ 
\end{defn}

\noindent One of the distinct advantages to these chain complexes is the ease with which we can prove theorems similar to those of E. S. Lee for the Khovanov homology of alternating links, \cite{Lee}. Not all surface tangle diagrams admit checkerboard colorings, but when they do we can repeat the definitions above for colored glyphs (see section \ref{sec:alternating}). This allows us to prove

\begin{theorem}[\ref{thm:lastone}] Let $T$ be a connected, alternating, checkerboard colored, tangle diagram for $(\Sigma, B)$. For each colored glyph $b^{c}$ with $b \in \glyph{\Sigma}{B}$ the chain complex $CT_{*}(T,b^{c})= \bigoplus_{i \in Z} CT_{i}(T,b^{c})$ is supported in a single homological grading when it is non-trivial. Consequently, $\underline{HT}_{\ast}(\mathcal{T}, b^{c}) \cong CT_{*}(T,b^{c})$ and any non-trivial homology groups occur in a single grading.
\end{theorem}

\noindent As a specific example, we recover a theorem similar to that in \cite{LR2} for the Khovanov homology of annular links:
\begin{cor}
Let $L$ be a connected, alternating link diagram in $\R^{2}\backslash\{(0,0)\}$ checkerboard colored so that the $1$-resolutions merge the black regions, and let $\mathcal{L}$ be the corresponding link in $\R^{3}$. Let $g_{k}$ be the glyph determined by the unit circle labeled with $k \in \Z$. Then any non-trivial groups $\underline{HT}_{\delta}(L,g_{k})$ is non-zero
only if
\begin{enumerate}
\item $k$ is odd, and $\delta + k = \sigma(\mathcal{L})$,
\item $k$ is even, the unbounded region is colored black, and $\delta + k = \sigma(\mathcal{L}) - 1$,
\item $k$ is even, the unbounded region is colored white, and $\delta + k = \sigma(\mathcal{L}) + 1$
\end{enumerate}
\end{cor} 
\noindent More examples can be found in section \ref{sec:alternating}.\\
\ \\
\noindent A natural question arises about the relationship between the homology theories of \APS in \cite{APS1}, \cite{APS2} and our own. First, note that \APS follow the conventions of \cite{Viro} while we follow those of \cite{Bar1}. Following T. Jaeger,\cite{TJ}, $\underline{HT}_{*}(\mathcal{T},b)$ will be identical to the $\delta$-graded (cohomology, due to the different convention) of those defined by \APS tensored with $\field{T}$ when $\mathcal{T}$ is a tangle without closed components.  \\
\ \\
\noindent Section \ref{sec:example} presents example calculations illustrating some of the phenomena that contribute to the boundary map. The reader should be able to understand the calculations after this introduction. Section \ref{sec:twisted} defines the totally twisted version of the complexes in \cite{APS1} and \cite{APS2}. We follow in section \ref{sec:collapse} with the proof that the totally twisted complexes are homotopy equivalent (over $\field{T}$) to the one described in this introduction. Section \ref{sec:Jaeger} proves that we may still employ the isomorphisms discovered by T. Jaeger, \cite{TJ}, allowing us to shift the location of the formal variables. We then use these isomorphisms to prove the invariance of the totally twisted chain complexes in section \ref{sec:invariance}. Section \ref{sec:alternating} proves the results about alternating projections references above. 

\section {Examples}\label{sec:example}

\diagram[0.5]{ex1}{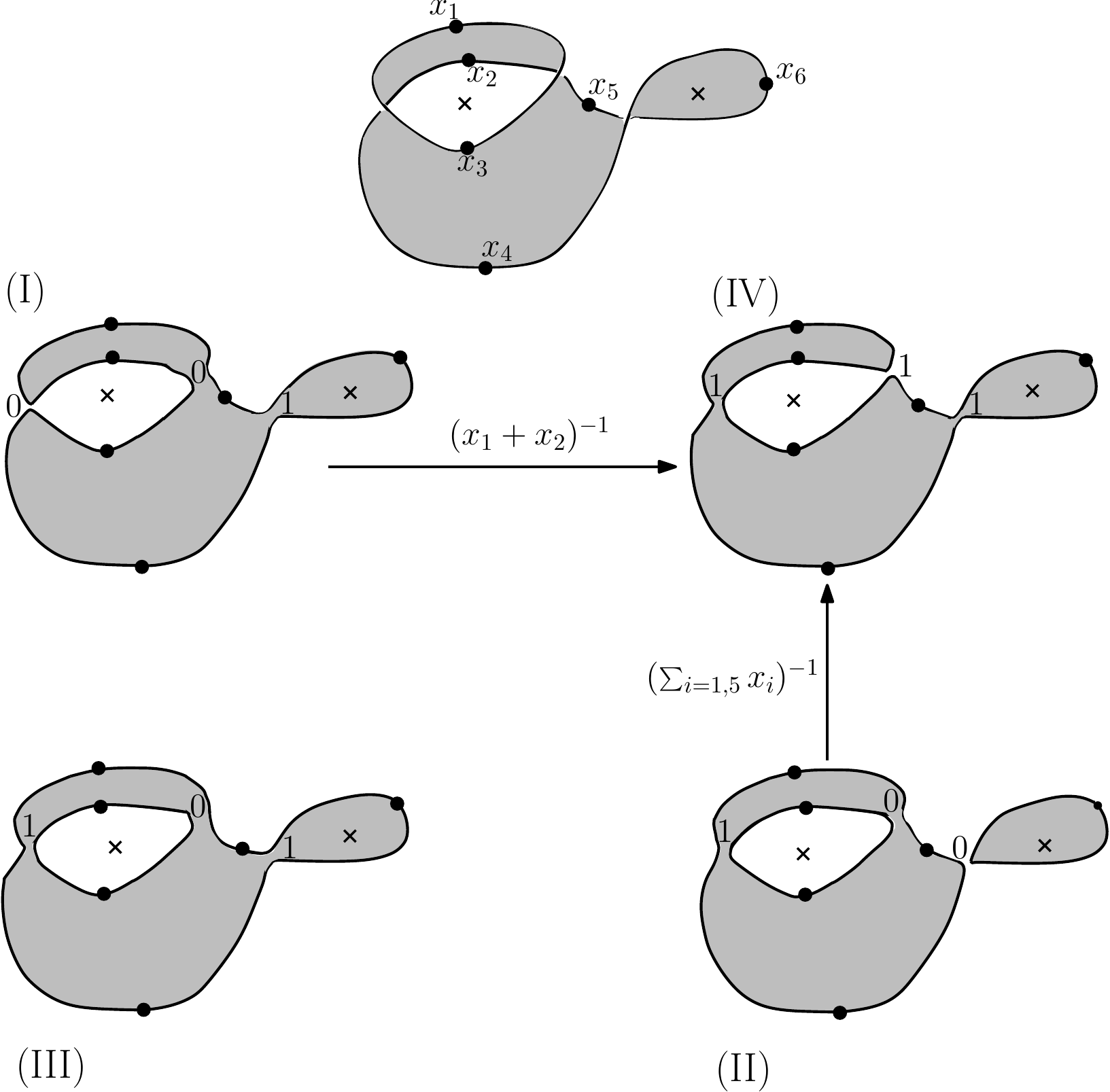}{An example of link diagram on a disk with two holes (see example 1).}
\diagram[0.5]{ex2}{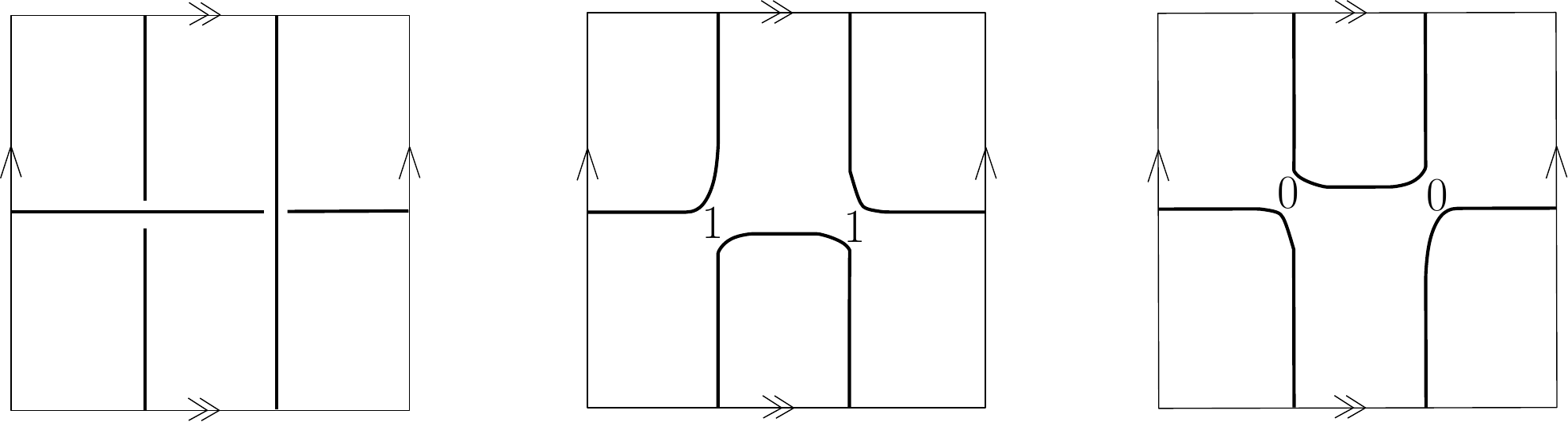}{An alternating link diagram on $T^{2}$ which cannot be checkerboard colored (see example 2).}
\diagram[0.5]{ex3}{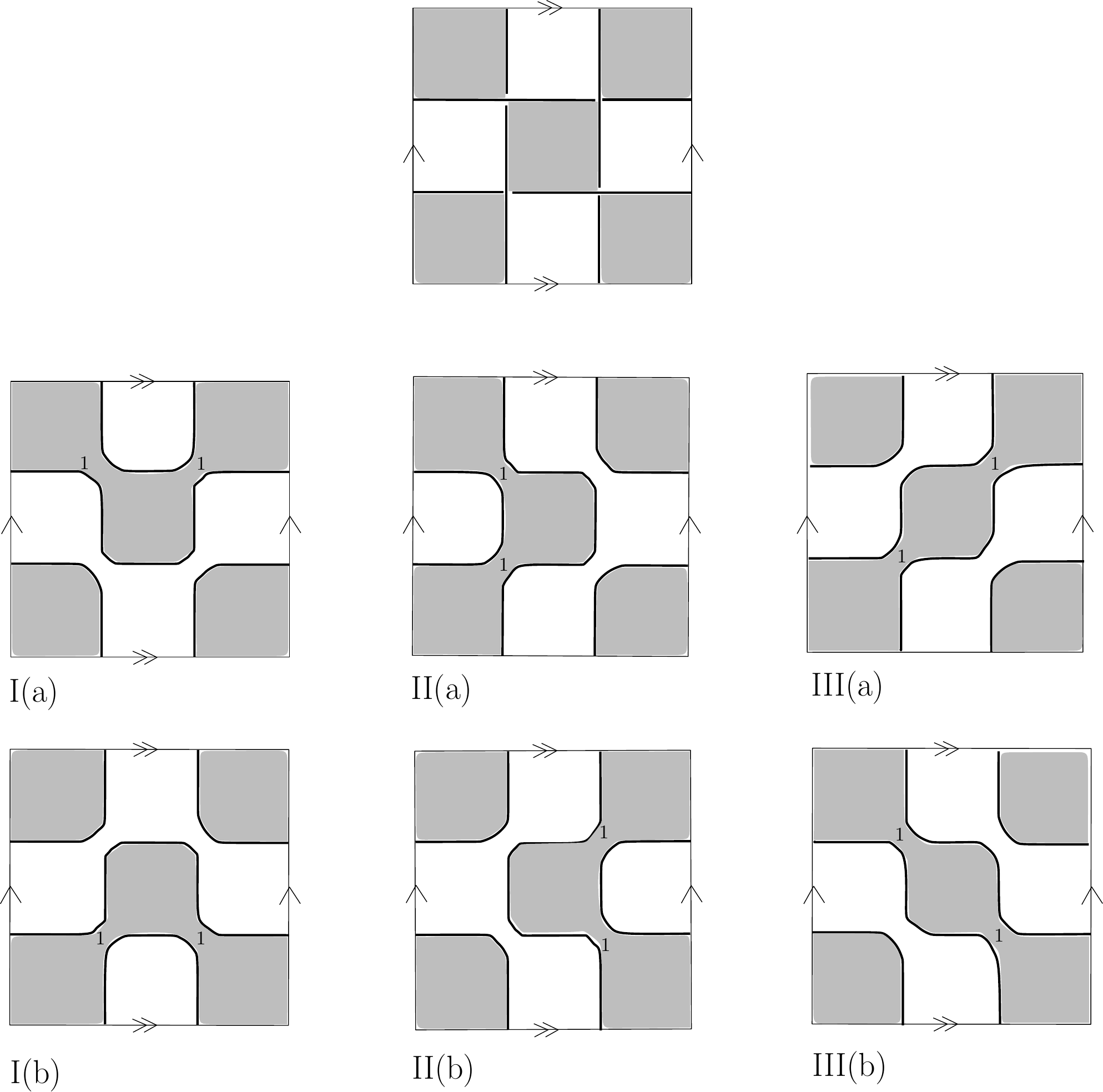}{An example of alternating link diagram on $T^2$. All the enhanced states without contractible circles have the same $k-$grading (see example 3).}
\diagram[0.5]{ex4}{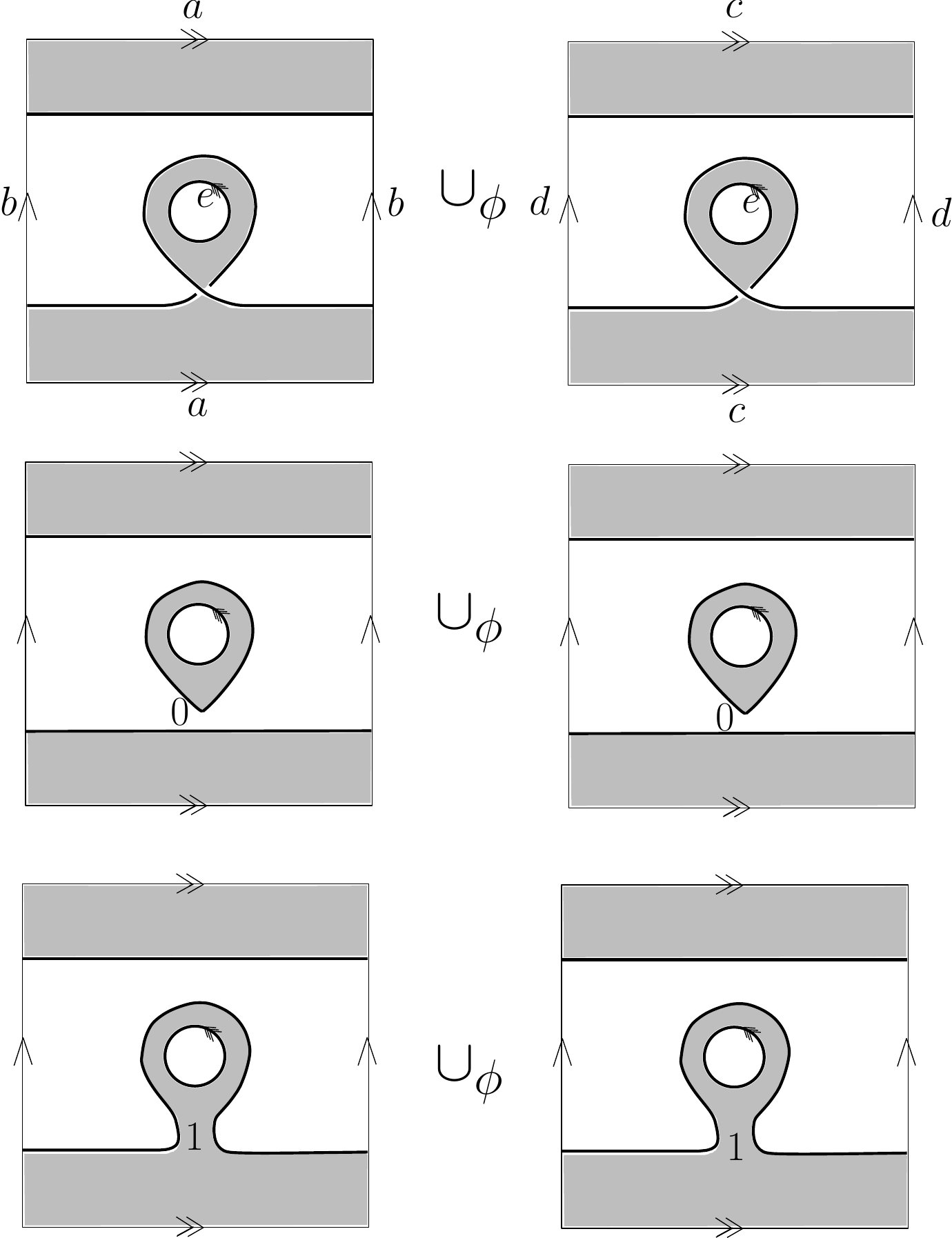}{An alternating link diagram on a closed genus $2$ surface where the coloring is necessary to collapse the chain complex to one grading. The two squares in each row have their vertical and horizontal sides glued in the normal manner for a torus. The holes labeled $e$ in the center of each square should be attached to each other by an identification map $\phi$ compatible with the orientations (see example 4).}
\diagram[0.75]{ex5}{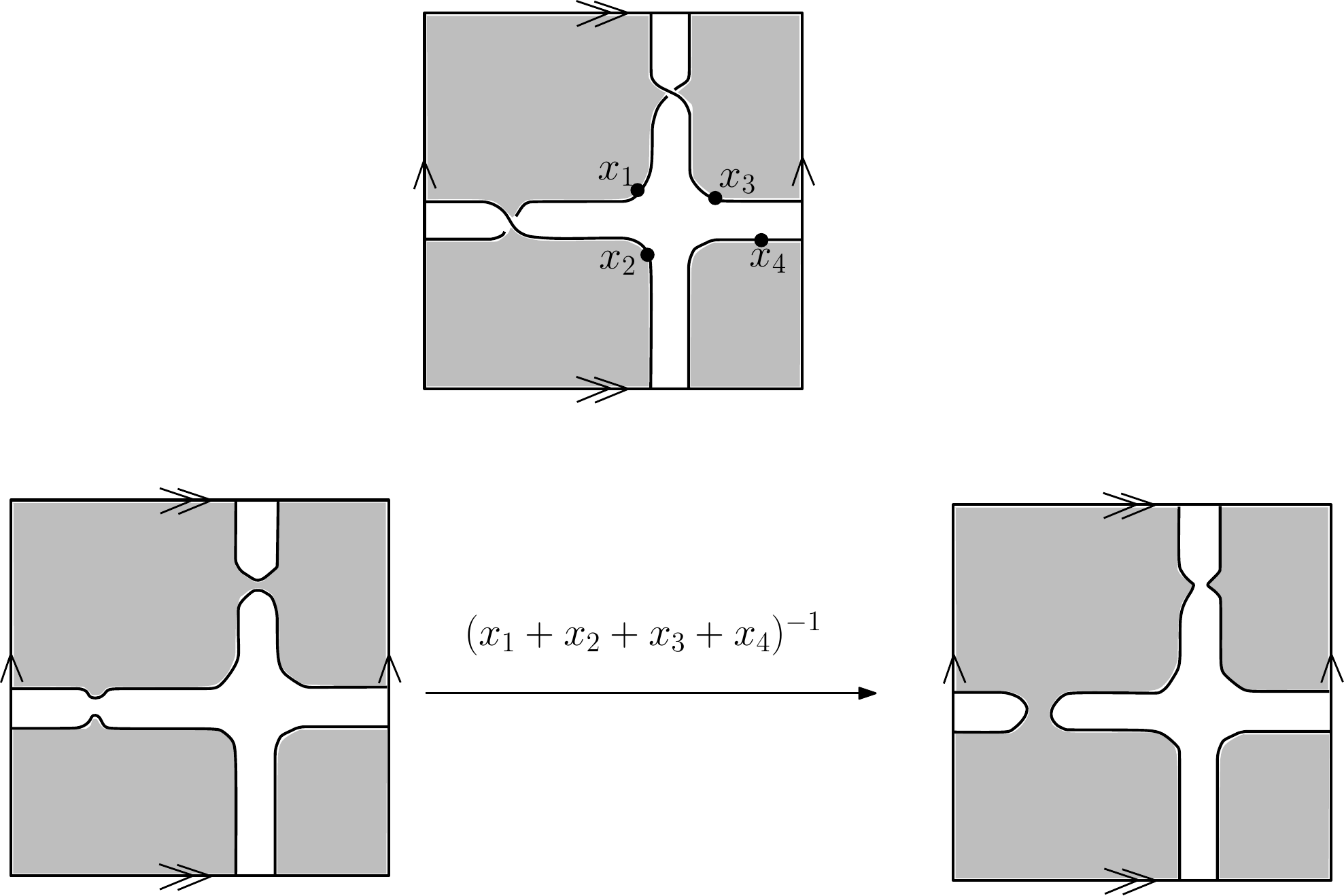}{A two component colored link diagram on a torus. The second row shows the only resolutions without contractible circles. If the pair of circles on the left are decorated with opposite signs, then there is a non-trivial boundary map to the resolution on the right. Note that the homotopy type of the noncontractible circles on the left is {\em different} from that on the right. This can only happen if for $0$-labeled circles in $\gly{r}$ (see example 5).} 

\noindent{\bf Example 1:}  Consider the diagram for a  two component link  depicted in the top diagram of Figure \ref{fig:ex1}, where the x's represent the points $p_{1}$ and $p_{2}$ (with $p_{1}$ on the left). We describe the complex for this diagram considered in $\R^{2}\backslash \{p_{1},p_{2}\}$, which provides a simple example of a non-trivial boundary map in the collapsed complex. \\
\ \\
\noindent The four lower diagrams are the only resolutions which do not have contractible circles. Diagrams (I) and (IV) have a single contractible circle we denote by $\gamma_{\{2\}}$. The smaller circle in (III) is $\gamma_{\{1\}}$ while the bigger is $\gamma_{\{1,2\}}$. Then $-2\gamma_{\{1\}} + \gamma_{\{2\}}$ is shorthand for the glyph assigning $-2$ to $\gamma_{\{1\}}$ and $+1$ to $\gamma_{\{2\}}$. Inspecting Figure \ref{fig:ex1} shows that there is only one state in (II) (assigning $-$'s to the left circles, and a $+$ to the right) which contributes to the chain complex for this glyph. Indeed all the glyphs in 
$\{\gamma_{\{1\}} \pm \gamma_{\{1,2\}}; -\gamma_{\{1\}} \pm \gamma_{\{1,2\}}; 2\gamma_{\{1\}} \pm \gamma_{\{2\}};-2\gamma_{\{1\}} \pm \gamma_{\{2\}}\}$ correspond to only one state among all those generated by these four diagrams.  As a result, $\widetilde{HT}_{*}(L, g)= \field{L}$ for glyphs $g$ in this set. The $\delta$-grading is readily computed for these generators: 
$$
\begin{array}{lll}
\gamma_{\{1\}} \pm \gamma_{\{1,2\}} & \lra & \delta=2 - (1\pm 1)  = 1 \mp 1 \\
-\gamma_{\{1\}} \pm \gamma_{\{1,2\}} & \lra & \delta=2 - (-1\pm 1)  = 3 \mp 1 \\
2\gamma_{\{1\}} \pm \gamma_{\{2\}} & \lra & \delta=1 - (2\pm 1)  =  -1 \mp 1 \\
-2\gamma_{\{1\}} \pm \gamma_{\{2\}}\} & \lra & \delta=1 - (-2\pm 1)  = 3 \mp 1 \\
\end{array}
$$
\ \\
However, when $g= \pm \gamma_{\{2\}}$ there are four states: one each from (I) and (IV) when we decorate the circle with $+$ or $-$, and two from (III) when we decorate the two $\gamma_{\{1\}}$ circles with different signs and place $+$ or $-$ on the $\gamma_{2}$ circle. Thus, the complexes $CT_{*}(L, \pm g)$ will be:
$$
0 \longrightarrow \field{L} \oplus \field{L} \oplus \field{L} \stackrel{\partial}{\longrightarrow} \field{L} \longrightarrow 0
$$ 
where $\partial(1,0,0) = (x_1+x_2)^{-1} \cdot v_{\pm}$, $\partial(0,1,0)=\partial(0,0,1)= (\sum_{i=1,5} x_i)^{-1} \cdot v_{\pm} $. Note that to compute $\partial(0,1,0)$ we needed to divide an annulus then join the resulting contractible circle to a noncontractible circle. Since $\field{L}$ is a field, $\partial$ is onto and $\mathrm{dim}(\mathrm{ker}\,\partial)=2$. Thus $\widetilde{HT}_{*}(L, \pm g)= \field{L}^{2}$ occurring in $\delta$-gradings $1-(\pm 1) = 1\mp 1$. \\
\  \\
\noindent{\bf Example 2:} Unlike for link diagrams in $\R^{2}$, many diagrams for links in $\Sigma \times (-1,1)$ cannot be checkerboard colored. Figure \ref{fig:ex2} depicts a simple example in $T^{2}$, where $T^2$ is presented by the standard identifications of the sides of a square. The $00$ and $11$-resolutions consist of a noncontractible circle, $\gamma$. Both resolutions contribute generators to the complex for the glyph $g$ where $\gamma$ is labeled with $+1$. The chain complex for this glyph $CT_*(L,g)$ is  $$0 \longrightarrow \field{L} \stackrel{0}  \longrightarrow \field{L} \longrightarrow 0$$ The middle map is $0$ since neither the $01$ nor $10$ resolutions have a contractible circle. Thus, there are
two generators for $\underline{HT}_{*}(L,g)$ with different homological gradings. Consequently, without the checkerboard coloring theorem \ref{thm:lastone} will not be true. \\
\ \\
\noindent{\bf Example 3:} The top diagram in Figure \ref{fig:ex3} show an alternating link diagram on $T^2$ with a given checkerboard coloring. Theorem \ref{thm:lastone} applies to this link, and we will verify its conclusion. The six diagrams
underneath the link diagram are the only resolutions without contractible circles. Each occurs with two $0$-resolved crossings and two $1$-resolved crossings, and thus all occur for the same homological grading. Since every state has the same homological grading, the boundary map is identically $0$ as predicted by Theorem \ref{thm:lastone}. The noncontractible circles in I(a) and I(b) are all homotopic, and will be denoted by $\gamma_{1}$. $\gamma_{2}$ will denote the circles in II(a) and II(b). However, III(a) and III(b) have distinct types of noncontractible circles,$\gamma_{3}$ and $\gamma_{4}$. It is straightforward to see that $\widetilde{HT}_{\ast}(L, \pm 2\gamma_{i})=\F^{2}_{(2 \mp 2)}$ for $i=1,2$ and $\widetilde{HT}_{\ast}(L, \pm 2\gamma_{i})=\F_{(2 \mp 2)}$ for $i=3,4$. On the other hand, each diagram contributes two generators to the trivial glyph, where one circle is assigned $+$ and the other is assigned $-$. Each of these generators occurs in $\delta$-grading $2(2) - 0 - 2 = 2$, Consequently, the homology for the trivial glyph is $\widetilde{HT}_{\ast}(L, \emptyset)=\F^{12}_{(2)}$.\\
\ \\
\noindent{\bf Example 4:} There are alternating diagrams in closed genus $n > 1$ surfaces where the glyph alone will not determine the homological grading, but the colored glyph does. Consider the alternating link diagram $L$ shown in Figure \ref{fig:ex4} on a closed genus $2$ surface, equipped with a checkerboard coloring. Note that a $1$-resolution at a crossing
will join the black quadrants. The all $0$ resolution is shown in the second row;  the noncontractible circles form three pairs, each cobounding an annulus. We consider the states where the circles in these pairs are decorated with different signs.
There are two white pair of pants, and three black annuli. To find the glyph we collapse the black annuli to obtain a white surface with three $0$ labeled circles. Thus the colored glyph is the trivial {\em white} glyph.  All $8$ of these states occur in  homological grading $0$. On the other hand, the all $1$ resolution, shown in the bottom row, has two pairs of noncontractible circles, cobounding white annuli, and a black region homeomorphic to a disc with four subdiscs removed. If we look at the four states where each pair of circles is decorated with different signs, we obtain generators corresponding to the trivial {\em black} glyph (found by collapsing the white annuli). These all occur in homological grading $2$. The Euler characteristic of the black regions in this case is $-2$. Note that in both cases the homological grading and the Euler characteristic of the black regions in the colored glyph add to $0$.\\
\ \\
\noindent{\bf Example 5:} Consider the diagram $L$ for a two component link in $T^{2} \times I$ shown at the top of Figure \ref{fig:ex5}. Calculating the boundary map for this example illustrates an important phenomenon. To start we identify the generators of the chain complexes. The only resolutions consisting of noncontractible circles are shown in the second row of Figure \ref{fig:ex5}. The diagram on the left is the all $0$ resolution, while the diagram on the right is the all $1$ resolution.   There are four generators in the chain complex for the trivial glyph $g_{\emptyset}$. Two are found by decorating the pair of noncontractible circles on the left with different signs, while the other two arise similarly from the diagram on the right. Thus, our chain complex $\widetilde{CT}_*(L,g_{\emptyset})$ will be: 
$$
0 \longrightarrow \field L \oplus \field L \stackrel{\partial}{\longrightarrow} \field L \oplus \field L \longrightarrow 0
$$
To compute $\partial$ we look at the resolutions with a single $1$-resolved crossing. The resolution where the top crossing is $1$-resolved and the left crossing is $0$-resolved has a single contractible circle $C$ through all the marked points. Thus $[C] = x_{1} + x_{2} + x_{3} + x_{4}$. Consequently, this choice of $1$-resolutions will contribute $(x_{1} + x_{2} + x_{3} + x_{4})^{-1}$ to the coefficient of {\em each} of the all $1$-resolution generators for $g_{\emptyset}$ in the boundary of each of the all $0$ resolution generators. Notice that the homotopy type of the noncontractible circle has also changed.  On the other hand, the resolution with a $1$-resolved crossing on the left and a $0$ resolved crossing at the top will, by the same reasoning, also contribute $(x_{1} + x_{2} + x_{3} + x_{4})^{-1}$ to the coefficient. Thus,
$$
\pairing{\partial (1,0)}{(1,0)} =  (x_{1} + x_{2} + x_{3} + x_{4})^{-1} + (x_{1} + x_{2} + x_{3} + x_{4})^{-1} = 0
$$
and $\partial \equiv 0$. Thus, $\widetilde{HT}_0(L, g_{\emptyset}) \cong \widetilde{HT}_2(L,g_{\emptyset}) \cong \field L \oplus \field L$. Note also that these will provide generators in different $\delta$-gradings in $\underline{HT}_{\ast}(L,g_{\emptyset})$. For completeness, we identify the homologies for the other glyphs, ignoring the $\delta$-gradings. Let $\gamma$ be a single horizontal circle and $\eta$ be a single vertical circle. Then  $\widetilde {HT}_*(L,g) \cong \field L$ for the glyphs $g = \pm 2 \gamma$ or $g = \pm 2 \eta$. All the other homologies are trivial.

\section{Defining the twisted skein complexes}\label{sec:twisted}

\subsection{Resolutions of $I$-bundle links and tangles}


\noindent Let $T$ be a (proper) diagram for an {\em oriented} tangle $\lk{T}$ in an I-bundle $\Sigma \times (-1,1)$ where $\Sigma$ is an oriented surface (see \cite{APS2} for more detail). The set of crossings in $T$ will be denoted $\cross{T}$, and the set of arcs will be denoted $\arcs{T}$. The number of positive crossings will be denoted $n_{+}(T)$ and the number of negative crossings will be $n_{-}(T)$. We will often omit the reference to $T$ when the choice of diagram is clear.\\
\ \\
\noindent As  in \cite{LR}, \cite{TJ}, we label each arc $f \in \arcs{T}$  with a formal variable $x_{f}$ and form the polynomial ring $\ring{T}=$ $Z_2[x_f|f\in \arcs{T}]$. The field of fractions of $\ring{L}$ will be denoted $\field{L}$.\\

\begin{defn}
For each subset $S \subset \cross{T}$ the {\em resolution} $r_{S}$ of $T$ is the surface tangle in $\Sigma$, considered up to isotopy, found by locally replacing each crossing $s \in \cross{T}$ according to the rule
$$
\inlinediag{L0} \hspace{.5in} \inlinediag{Linfty} 
$$ 
The set of resolutions for $T$ will be denoted $\res{T}$. For each resolution $r$, denote by $|r|$ the number of elements in the corresponding subset $S \subset \cross{T}$. 
\end{defn}

\noindent{\bf Note:} Given a crossing $c \in \cross{T}$ and a resolution $r=r_{S}$ we will also use the notation $r(c) = 0$ for $c \not\in S$ and $r(c)=1$ for $c \in S$. Thus $r$ will stand both for the resolution diagram and for the indicator function for the set of crossings defining the resolution.\\
\ \\
\noindent We assign a weight to each circle (or arc) in a resolution by adding the formal variables along the arcs composing the circle (or arc). However, we will only use the weights on the contractible circles. \\

\begin{defn}
Let $r$ be a resolution for $T$. By $\circles{r}$ we denote the set of circles in the resolution $r$. For $C \in \circles{r}$, define
$$
[C] = \begin{cases} 
 \sum_{f \in \arcs{C}} x_{f}& \text{if $C$ bounds a disc in $\Sigma$}\\
 0& \text{if $C$ does not bound a disc in $\Sigma$}
\end{cases}
$$
Circles which bound a disc in $\Sigma$ are contractible, and will be referred to as {\em contractible circles}. The other circles will be called {\em noncontractible}. The set of contractible circles will be denoted $\circles[\bullet]{r}$, while the noncontractible circles will be denoted $\circles[\circ]{r}$.
\end{defn}

\subsection{The Koszul complex of a resolution} For each resolution $r$ of a tangle diagram $T$ we associate a chain complex which will be called the vertical complex of $r$. We do this in several stages. For more detail see the introduction of  \cite{LR}.\\
\ \\
\noindent To a circle $C \in \circles{T}$ we associate the complex $\mathcal{K}([C])$ defined as
$$
0 \longrightarrow \field{T}\,v_{+} \stackrel{\cdot [C]}{\longrightarrow} \field{T}\,v_{-} \longrightarrow 0
$$
where $v_{\pm}$ occur in gradings $\pm 1$. The differential in $\mathcal{K}([C])$ will be denoted $\partial_{C}$. \\
\ \\
\noindent We use these complexes to construct a Koszul complex for the resolution $r$. Let $\circles{r}=\{C_1,...,C_k\}$, then 
$$
\mathcal{K}(r) = \mathcal{K}([C_{1}], \ldots, [C_{k}]) = \mathcal{K}([C_{1}]) \otimes_{\field{T}} \mathcal{K}([C_{2}]) \otimes_{\field{T}} \cdots \otimes_{\field{T}} \mathcal{K}([C_{k}])
$$
The differential in $\mathcal{K}(r)$ is then $\partial_{\mathcal{K}(r)} = \sum_{i > 0} \partial_{C_{i}}$. Since we only define the complex in characteristic 2, the ordering of the factors does not matter.\\
\ \\
\noindent In the final step, we will shift the gradings in $\mathcal{K}(r)$.

\begin{defn}
Let $M = \oplus_{\vec{v} \in \Z^{k}} M_{\vec{v}}$ be a $\Z^{k}$-graded $R$-module, then $M\{\vec{w}\}$ is the $\Z^{k}$-graded module with $(M\{\vec{w}\})_{\vec{v}} \cong M_{\vec{v}-\vec{w}}$. 
\end{defn}
\noindent Consequently, $\mathrm{grad}(m\{\vec{w}\}) = \mathrm{grad}(m) + \vec{w}$ on homogeneous elements $m \in M$.\\

\begin{defn}
The vertical complex for the resolution $r$ of the tangle diagram $T$ is the complex
$$\mathcal{V}(r) = \mathcal{K}(r)\{|r|\}$$
\end{defn}

\noindent The grading on $\mathcal{V}(r)$ for any resolution  will be called the quantum grading. The differential $\partial_{\mathcal{V}(r)}$ changes the quantum grading by $-2$.\\
\ \\
\noindent As in \cite{APS1}, we represent the basis elements for $\mathcal{V}(r)$ by decorating the diagram $r$: a basis element corresponds to a choice of $+$ or $-$ for each circle $C_{i}$, $i > 0$. The basis element is then the tensor product of $v_{+}$ and $v_{-}$ elements where the $i^{th}$ factors is determined by the sign on $C_{i}$, \cite{Bar1}. We will call these basis elements, and their corresponding diagrams, {\em pure states}, while elements of $\mathcal{V}(r)$ will be called states.\\ 

\begin{defn} 
\noindent For each pure state of $r$ we define
\begin{enumerate}
\item the {\em homology grading} to be $h(r) = |r|$,
\item the {\em quantum grading} to be $q(r) = |r| + \#(+ \text{ circles}) - \#(-\text{ circles})$, and
\item the {\em surface grading} to be $k(r) = \#(+ \text{ noncontractible circles}) - \#(- \text{ noncontractible circles})$
\end{enumerate}
\end{defn} 

\noindent For the bi-grading $r \rightarrow (h(r), q(r))$ the map $\partial_{\mathcal{V}(r)}$ is a $(0,-2)$ differential. 

\subsection{Adding the horizontal differential}

\noindent Following the construction in \cite{APS1}, \cite{APS2}, we combine the complexes $\mathcal{V}(r)$ into a triply graded module $\widetilde{SK}(T)$ by setting  
$$
K_{i}(T) = \bigoplus_{r \subset \res{T},\ |r| = i}  \mathcal{V}(r) \hspace{0.5in} i \in \Z
$$
and taking the direct sum
 $$\widetilde{SK}(T) = \bigoplus_{i \in \Z} K_{i}(T)$$
$\widetilde{SK}(T)$ is a $q$-graded chain complex with differential 
$$\partial_{\mathcal{V}}= \bigoplus_{r \in \res{T}} \partial_{\mathcal{V}(r)}$$

\noindent The other gradings on $\widetilde{SK}(T)$ are given by the homology grading, $K_{i}(T) \rightarrow i$, and the surface grading. The homogeneous elements of $\widetilde{SK}(T)$ in homology grading $i$, quantum grading $j$, and surface grading $k$ will be denoted $\widetilde{K}^{(i,j,k)}$.  Shifts in the triple grading will be denoted by $\{(\Delta i, \Delta j, \Delta k)\}$. \\
\ \\ 
\noindent \APS in \cite{APS1}, \cite{APS2}, discovered a $(+1,0,0)$ differential, $\partial_{\Sigma}$, on this triply graded module, which can be described in terms of the pure states by figures \ref{fig:Khovanov-map1}, \ref{fig:Khovanov-map2}, and \ref{fig:Khovanov-map3}.
\diagram[0.55]{Khovanov-map1}{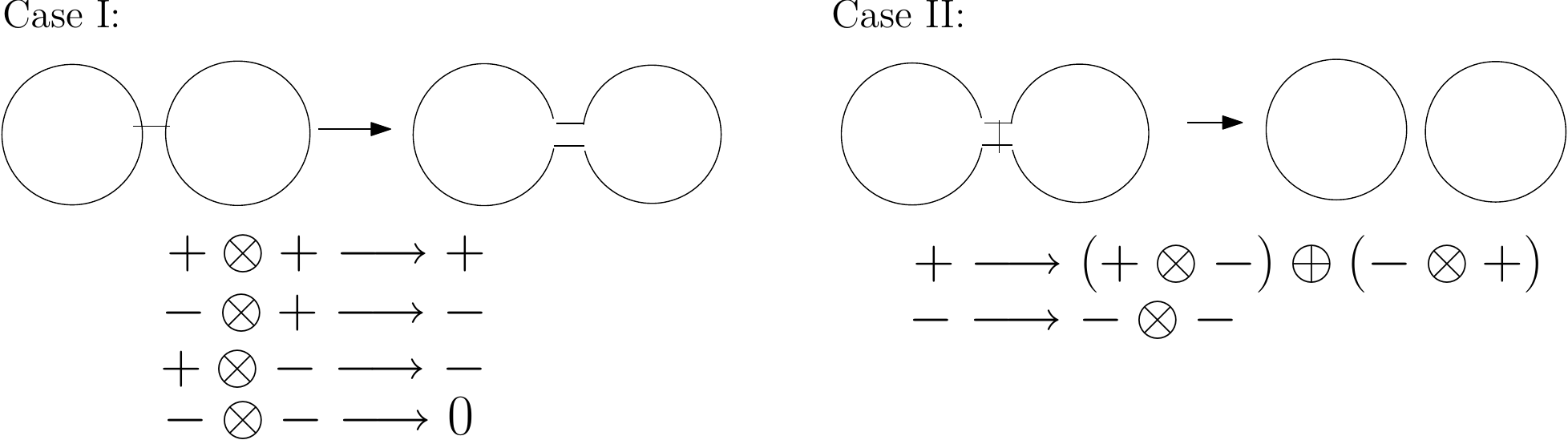}{Cases I and II are the standard maps from Khovanov homology, \cite{Khov}. Each row indicates a non-zero entry in $\partial_{\Sigma}$. Given a pure state $r_{S}$ with components as on the left, and the indicated crossing $c \not\in S$, there is an non-zero term $r_{S'}$ in $\partial_{\Sigma}$
where $r_{S'}$ is the pure state on the right. The first (second) factor of the tensor product corresponds to the left (right) circles in the diagram.} 

\diagram[0.55]{Khovanov-map2}{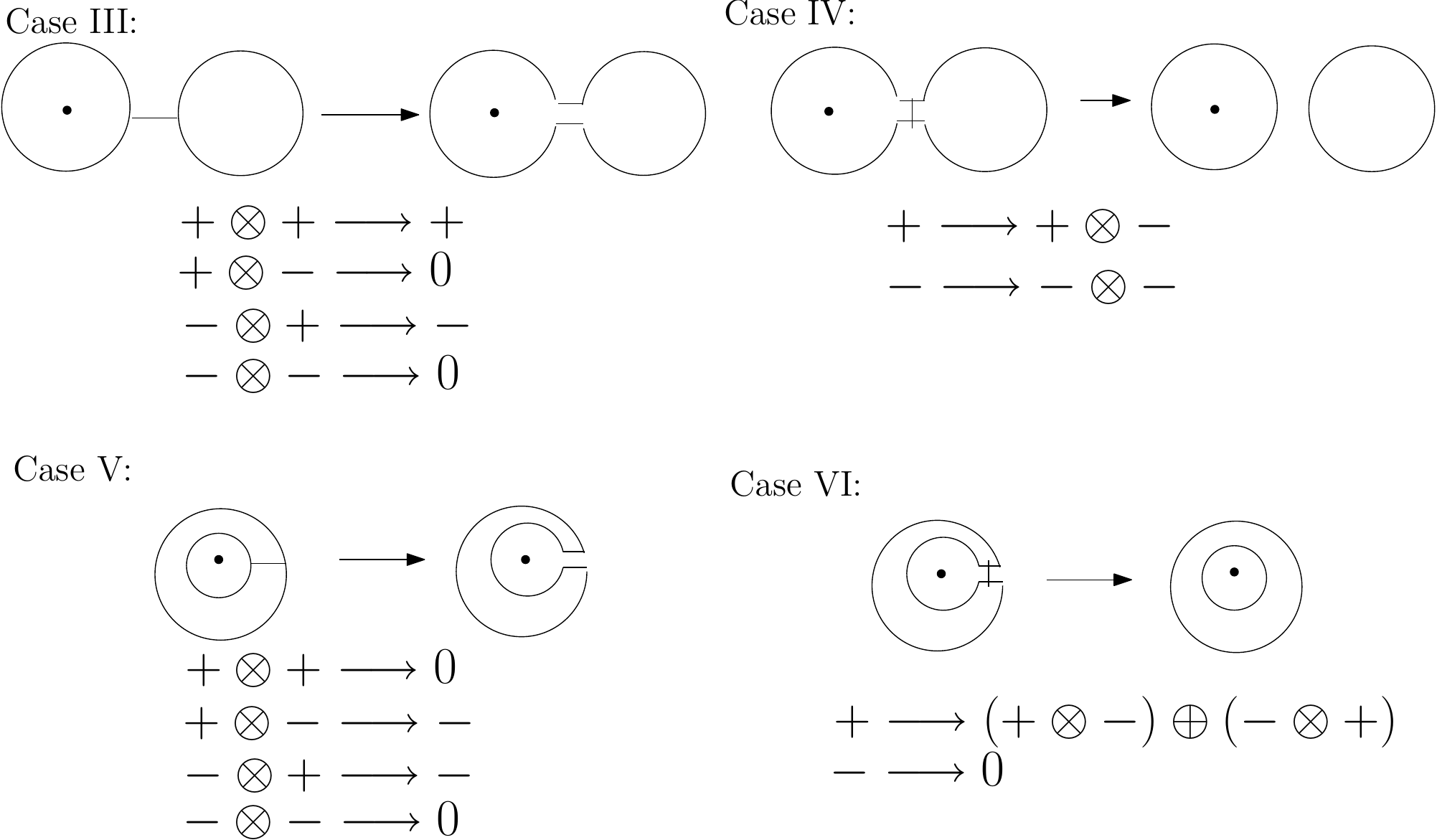}{Each row indicates a non-zero entry in $\partial_{\Sigma}$ for a noncontractible circle, as in [APS1]. noncontractible circles are represented by circles surrounding a dot, while contractible circles contain no dot. Two noncontractible circles surrounding the same dot should be homotopic to each other in $\Sigma$; and conversely, if they have different dots, they are not homotopic. For nested circles, the first factor in the tensor product corresponds to the outside circle.} 

\diagram[0.55]{Khovanov-map3}{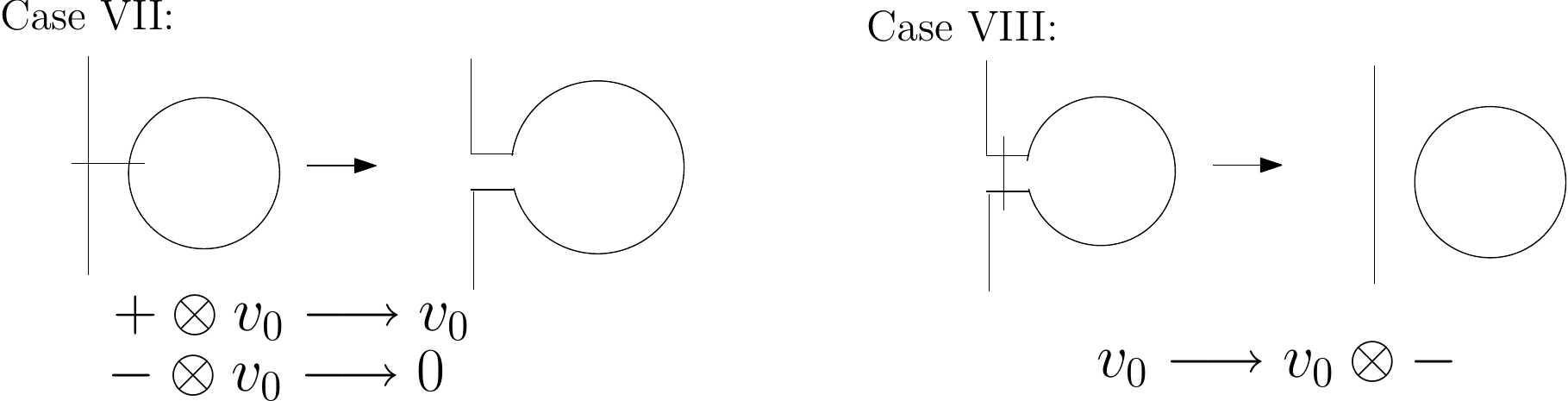}{Each row indicates a non-zero entry in $\partial_{\Sigma}$ as in [APS2]. The left segment corresponds to an arc in the planar resolution of a tangle. The circle is {\em always} contractible. All arc are decorated with  a special element $v_{0}$.} 

\noindent Our next step is to show that $\partial_{\Sigma}$ commutes with $\partial_{\mathcal{V}}$ on $\widetilde{SK}(T)$.

\begin{prop} \label{prop:complex}
The maps $\partial_{\Sigma}: \widetilde {SK}^{(i,j,k)}(T) \longrightarrow  \widetilde {SK}^{(i+1,j,k)}(T)$ and
$\partial_{\mathcal{V}}: \widetilde {SK}^{(i,j,k)}(T) \longrightarrow  \widetilde {SK}^{(i,j-2,k)}(T)$ satisfy
$$
\partial_{\Sigma}\,\partial_{\mathcal{V}} + \partial_{\mathcal{V}}\,\partial_{\Sigma} = 0
$$
\end{prop}

\noindent The proof of proposition \ref{prop:complex} will be deferred to the next section. However, proposition \ref{prop:complex} ensures that we can consider $\widetilde{SK}(T)$ as a double chain complex where with the horizontal differential being $\partial_{\Sigma}$, the differential discovered by \APS, and with the vertical differential being $\partial_{\mathcal{V}}$, the Koszul differential. By collapsing the first two gradings using $\delta(i,j) = 2i - j$ we can make $\partial_{\mathcal{V}} + \partial_{\Sigma}$ into a $(+2,0)$ differential. \\

\begin{defn}
The {\em totally twisted skein homology complex} for a tangle digram $T$ in $\Sigma$ is the $\field{T}$ vector space
$\underline{SK}(T):=\widetilde{SK}(L)\{(-n_{+}(L),0,0)] $
equipped the bigrading $(\delta, s): \underline{SK}^{(i,j,k)}(T) \rightarrow (2i - j, k)$ and the $(+2,0)$ differential $\underline{\partial} = \partial_{\mathcal{V}} + \partial_{\Sigma}$
\end{defn}

\subsection{Refinements}

\noindent Just as in \cite{APS2}, the complex $\underline{SK}(T)$ is a direct sum of subcomplexes in several ways. First, since both $\partial_{\Sigma}$ and $\partial_{\mathcal{V}}$ preserve the surface grading the complex splits as
$$
\underline{SK}(T) = \bigoplus_{k \in \Z} \underline{SK}(T,k)
$$
where $\underline{SK}(T,k) = \bigoplus_{\delta \in \Z} \underline{SK}_{\delta, k}$ equipped with the restriction of $\underline{\partial}$ as the differential. \\
\ \\
\noindent Actually, as in \cite{APS1}, we can refine the surface grading further. For this we need some additional data 
\begin{defn}
A {\em glyph} for  $(\Sigma, B)$ is an equivalence class of disjoint proper embeddings in $\Sigma$ of
\begin{enumerate}
\item a finite collection of noncontractible circles $C_{i}$, each labeled with an integer $n_{i} \in \Z$, and
\item a finite collection of arcs, whose boundaries exhaust $B$   
\end{enumerate}
considered up to the equivalence generated by boundary preserving surface isotopy, and the erasure/introduction of a $0$ labeled noncontractible circles. We will denote the set of glyphs by $\glyph{\Sigma}{B}$. In $\glyph{\Sigma}{\emptyset}$ there is a {\em trivial glyph}, the class where the number of circles embedded is $0$. 
\end{defn}
\noindent{\bf Note:} The trivial glyph is also the class represented by disjoint embeddings of noncontractible circles in $(\Sigma, \emptyset)$ all of whom are labeled with $0$.\\
\ \\
\begin{defn}
For a glyph $g$ in $\glyph{\Sigma}{B}$, $k(g)$ is the sum of the labels of all the noncontractible circles in $g$.
\end{defn}

\begin{defn}
For each pure state $r$ for a tangle diagram $T$, $\gly(r) \in \glyph{\Sigma}{B}$ is the glyph found by first erasing all the contractible circles in $r(T)$, and replacing each collection of noncontractible circles that pairwise cobound annuli by a single circle whose label is the number of $+$ circles in the set minus the number of $-$ circles. 
\end{defn}
\ \\
\noindent From Figures \ref{fig:Khovanov-map1}, \ref{fig:Khovanov-map2}, \ref{fig:Khovanov-map3}  we observe that the only non-trivial entries in $\partial_{\Sigma}$ either 1) cleave off or merge a contractible circle while leaving the decorations on arcs and noncontractible circles alone, or 2) divide/introduce a $+,-$ pair of noncontractible circles which cobound an annulus. Consequently, if $\beta(r) = b$ and $\pairing{\underline{\partial} r}{r'} \neq 0$, then $\beta(r') = b$ as well. \\
\ \\
\noindent Let $\underline{SK}(T,b)$ be the sub-complex of $\underline{SK}(T)$ spanned by pure states with surface grading equal to $b$, equipped with the restriction of $\underline{\partial}$ and graded by $\delta$. Then 
$$
\underline{SK}(T) = \bigoplus_{b \in \glyph{\Sigma}{B}} \underline{SK}(T,b)
$$
as chain complexes.

\noindent In the section \ref{sec:invariance} we prove\\

\begin{theorem}\label{theorem:main}
For each $b \in \glyph{\Sigma}{B}$ the homology $H_{\ast}(\underline{SK}(T, b), \underline{\partial})$, as a $\delta$-graded module, is invariant under the first, second and third Reidemeister moves applied to $T$.   
\end{theorem}

\section{Proof of proposition \ref{prop:complex}}

\begin{prop}
Let $r_{S}$ be a resolution of $T$ and let $c \in \cross{T}$ with $r(c) = 0$. If $C$ is a circle in  $r_{S \cup \{c\}}$ formed by merging the circles $C_{1}$ and $C_{2}$ in $r_{S}$, then the map
in Figures \ref{fig:Khovanov-map1}, \ref{fig:Khovanov-map2} induces a chain map on the Koszul complexes (over $\field{T}$)
$$
\mu : \mathcal{K}(\area{C_{1}}) \otimes \mathcal{K}(\area{C_{2}})\longrightarrow \mathcal{K}(\area{C})\{1\}
$$
Likewise, if $C_{1}$ and $C_{2}$ arise in $r_{S \cup \{c\}}$ from dividing a circle $C$ in $r_{S}$, then
$$
\Delta : \mathcal{K}(\area{C}) \longrightarrow \big(\mathcal{K}(\area{C_{1}}) \otimes \mathcal{K}(\area{C_{2}})\big)\{1\}
$$
is a chain map over $\field{T}$. 
\end{prop}
\ \\
\noindent{\bf Proof:} This proposition for Cases I and II was addressed in \cite{LR} (see also \cite{TJ}). We just need to prove the result for cases III, IV, V, and VI.  Figure \ref{fig:pic2} illustrates the proof. By examining the diagrams we can see that in each of these four cases  $\partial_{\lk KH}\circ\partial_V=\partial_V \circ \partial_{\lk KH}=0$ even though $[C]=0$ if $C$ is a noncontractible circle. Note that several of the horizontal arrows are zero maps.$\Diamond$

\diagram[0.75]{pic2}{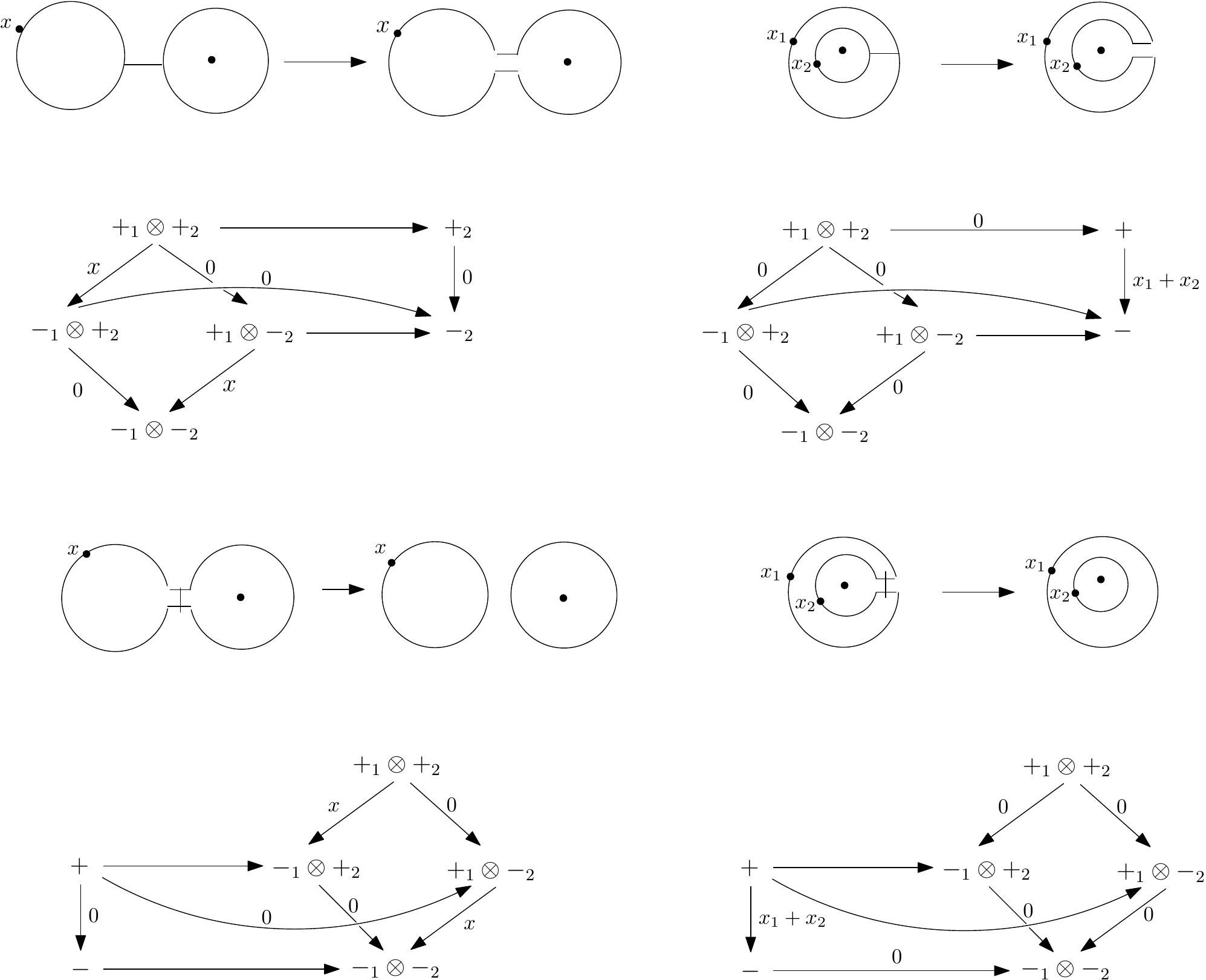}{} 

\diagram[0.75]{pic2b}{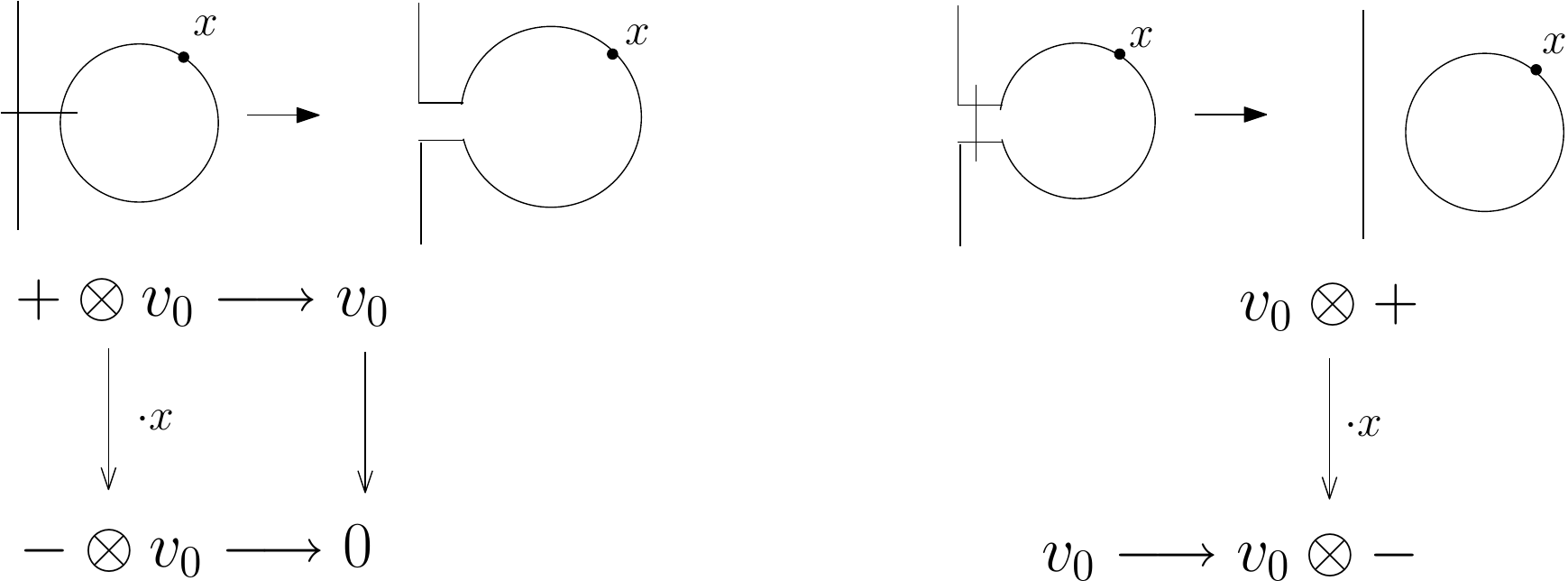}{}

\begin{prop}
Let $r_{S}$ be a resolution of $T$ and let $c \in \cross{T}$ with $r(c) = 0$. If $A'$ is a arc in  $r_{S \cup \{c\}}$ formed by merging a circle $C$ with an arc $A$ in $r_{S}$, then the map
in case VII of Figure \ref{fig:Khovanov-map3} induces a chain map on the Koszul complexes (over $\field{F}$)
$$
\mu : \field{T}v_{0}\otimes \mathcal{K}[C]  \longrightarrow \field{T}v_{0}\{1\}
$$
Likewise, if the circle $C$ and arc $A$ arise in $r_{S \cup \{c\}}$ from dividing an arc $A'$ in $r_{S}$ (as in case VIII), then
$$
\Delta : \field{T}v_{0} \longrightarrow  \field{T}v_{0} \otimes  \mathcal{K}[C]\{1\}
$$
is a chain map over $\field{T}$. 
\end{prop}

\noindent {\bf Proof:} For those cases involving an arc in the tangle $T$ (cases VII and VIII in Figure \ref{fig:Khovanov-map3}) 
we verify the property using Figure \ref{fig:pic2b}
\noindent $\Diamond$\\
\ \\
\noindent It is then simply a matter of taking tensor products to conclude that these maps induce chain maps on $\mathcal{K}(r)$ for each resolution $r$.\

\section{Collapsing the $\delta$-graded chain complex}\label{sec:collapse}

\noindent We have seen that $\partial =\partial_{\lk{V}}+ \partial_{\Sigma}$ is a differential which preserves the glyph of a state. Thus, we can fix the value of $b \in \glyph{\Sigma}{B}$ and let $\lk{O}(T,b) $  be the set of pure states whose glyph equals $b$. We let $\lk{O}_{n}(T,b)$ be the subset of $\lk{O}(T,b)$ consisting of those states that do not have any contractible circles. In this section we will describe a deformation retraction of $\widetilde{SK}(T,b)$  supported on $\lk{O}_{n}(T,b)$. We denote it $(CT_{*}(T,b), \partial_{n})$. \\
\ \\
\noindent Notice that each state $S$ in $\widetilde{SK}(T,b)$ with a contractible circle $C$ will have a non-trivial vertical differential arising from the comples $0 \lra \field{T}\,v_{+} \stackrel{[C]}{\lra} \field{T}v_{-} \lra 0$. Since the middle map is an isomorphism, we can cancel this portion of the differential (see \cite{LR} or any reference on computational homology) to simplify the complex by a chain homotopy equivalence removing any state with $C$ in its resolution. By repeatedly using this cancellation for pure states $S \in (\lk{O}(T,b) \backslash \lk{O}_n(T,b))$ we obtain a chain homotopy equivalent complex supported on $\lk{O}_n(T,b)$. Notice that since we are canceling along contractible circles, the glyph $b$ is unchanged. This procedure changes the boundary map, and we now proceed to analyzing that alteration. We presume the reader is familiar with cancellation in chain complexes, which is covered in many references, although its use here parallels the description in \cite{LR}.\\   
\ \\
\noindent Let $r$, $r' \in \lk{O}_n^b$ and suppose $\pairing{\partial_n(r)}{r'} \neq 0$. Then there is a sequence
of states $r=r^{+}_{0},r^{-}_{1},r^{+}_{1}, \ldots, r^{+}_{k-1}, r^{-}_{k} = r'$ where each transition $r^{+}_{i} \rightarrow r^{-}_{i+1}$ corresponds to a term in $\partial_{\Sigma}$ and each transition $r^{-}_{i+1} \rightarrow r^{+}_{i+1}$ comes from a Koszul map.\\
\ \\
\noindent Throughout this path of transitions we record the vector $(m_{+}(\overline{r}),m_{-}(\overline{r}))$ where $m_{\pm}$ is the number of positively/negatively decorated contractible circles in the state $r$. Thus we have $(m_{+}(r),m_{-}(r)) = (0,0) = (m_{+}(r'),m_{-}(r'))$. We also record the sequence of homological gradings for these states. The terms in $\partial_{\Sigma}$ and $\partial_{\mathcal{V}}$ alter these data by
\begin {enumerate}
\item A transition $r^{+}_{i} \rightarrow r^{-}_{i+1}$ changes the quantum grading by $-1$ and thus corresponds to changing
$(m_{+}, m_{-})$ by either $(-1,0)$ or $(0,1)$. It also changes the homological grading by $+1$.
\item A transition $r^{-}_{i+1} \rightarrow r^{+}_{i+1}$ changes the quantum grading by $+2$, the homological grading by $0$ and
the vector $(m_{+}, m_{-})$ by $(+1,-1)$ as we switch a $-$ decorated circle to a $+$.
\end{enumerate}
The sequence stops if we ever have $m_{-}(r^{-}_{i}) = 0$ since there will be no contractible circle whose decoration can be changed from $-$ to $+$. For the sequence above to exist we therefore need $m_{-}(r^{-}_{i}) > 0$ for $i=1, \ldots, k-1$. Assume for the moment that $k > 2$. Therefore, $(m_{+}(r^{-}_{1}), m_{-}(r^{-}_{1})) = (0,1)$ and $(m_{+}(r^{+}_{1}), m_{-}(r^{+}_{1})) = (1,0)$. We will show by induction that $(m_{+}(r^{+}_{i}), m_{-}(r^{+}_{i})) = (i,0)$ for $i = 2, \ldots, k-1$. Assume this is true for $i-1$. Then $(m_{+}(r^{-}_{i}), m_{-}(r^{-}_{i})) = (i-1,1)$ is determined by the requirement that $m_{-} > 0$. Therefore, $(m_{+}(r^{+}_{i}), m_{-}(r^{+}_{i})) = (i,0)$. For $i = k-1$ we then have $(m_{+}(r^{+}_{k-1}), m_{-}(r^{+}_{k-1})) = (k-1,0)$ which can only change by $(-1,0)$ or $(0,1)$ in the transition $r^{+}_{k-1} \rightarrow r^{-}_{k} = r'$. Thus
$(m_{+}(r'),m_{-}(r'))$ is either $(k-2,0)$ or $(k-1,1)$. Neither of these can be $(0,0)$ when $k > 2$. Thus, for $(m_{+}(r'),m_{-}(r')) = (0,0)$, we need to have $k = 2$. \\
\ \\
\noindent Furthermore, this argument gives a picture of how $\pairing{\partial_n(r)}{r'} \neq 0$ can arise. We must 
start with noncontractible circles and arcs only, change a single crossing from a $0$ resolution to a $1$ resolution while  producing a single contractible circle $C$ decorated with a $-$, then change the minus to a plus at the expense of multiplication by $1/[C]$, and then change another crossing from $0$ to $1$ which merges the $+$-decorated contractible circle into a noncontractible circle or arc, or splits $C$ into two noncontractible circles which cobound an embedded annulus. \\
\ \\
\noindent Thus, the new differential map $\partial_n$ changes the grading by $2$ and includes terms from each of the configurations in Figure \ref{fig:pic13}

\diagram[0.75]{pic13}{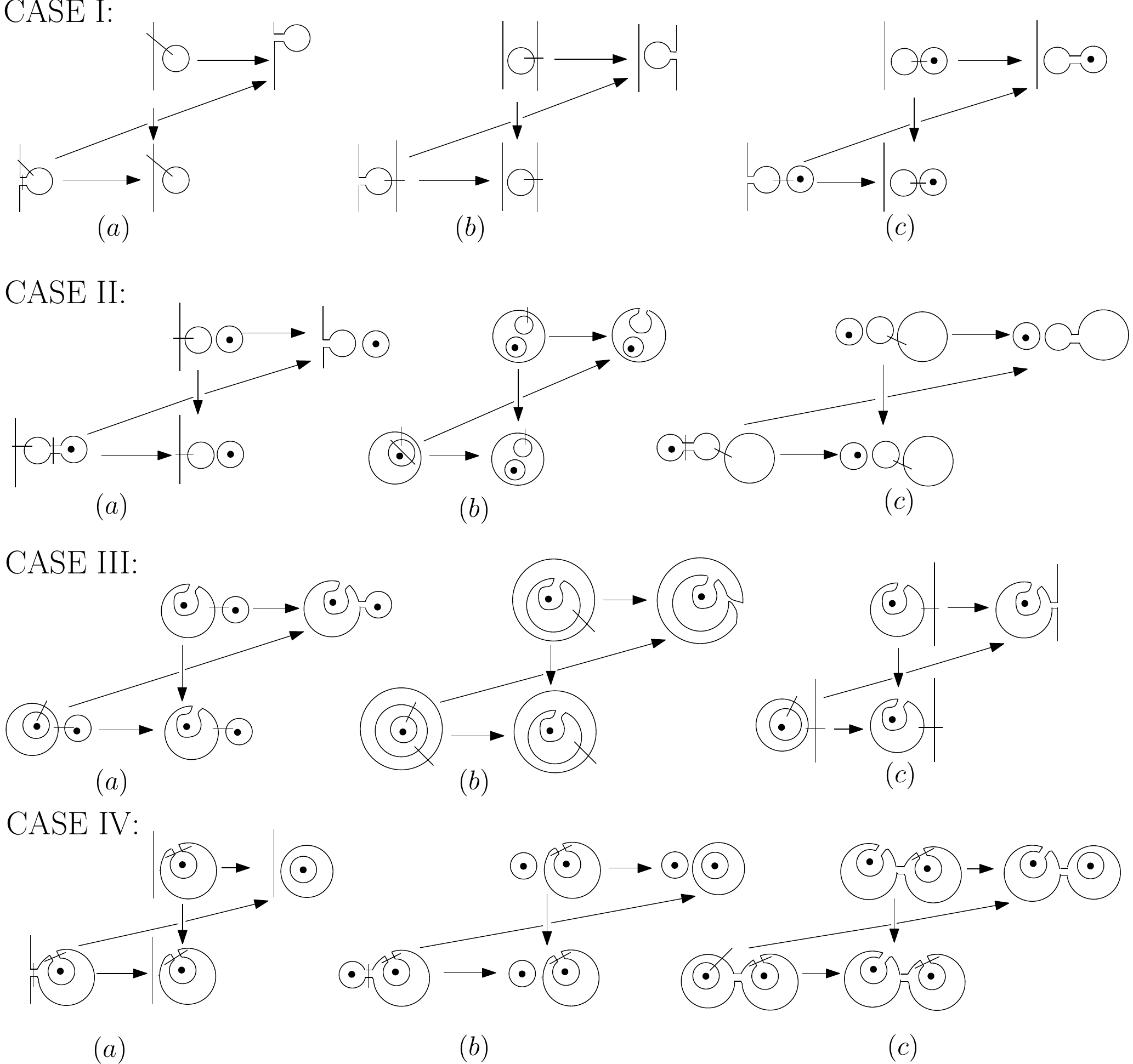}{} 

\begin{defn}
Suppose $r \in  \lk{O}_n(T,b)$ and $c_{i} \in \cross{T}, i =1,2$ with $r(c_{i}) = 0$. Let $r_{ab}$, $a,b \in \{0,1\}$ be the resolution with $r_{ab}(c_{1}) = a$, $r_{ab}(c_{2}) = b$ and $r_{ab}(c) = r(c)$ for every other crossing. When $\gly{r_{ab}}=b$ for $a,b \in \{0,1\}$, let
\begin{equation}
\pairing{r}{r_{11}} = \left\{\begin{array}{cl} 
\displaystyle{\frac{1}{[C_{10}]} + \frac{1}{[C_{01}]}} & r_{10}\mathrm{\ and\ }r_{01}\mathrm{\ have\ contractible\ circles\ }C_{10}\mathrm{\ and\ }C_{01}\mathrm{\ respectively} \\
\displaystyle{\frac{1}{[C_{10}]}} & r_{10}\mathrm{\ has\ a\ contractible\ circle\ }C_{10}\mathrm{\ but\ }r_{01}\mathrm{\ does\ not}\\
\displaystyle{\frac{1}{[C_{10}]}} & r_{01}\mathrm{\ has\ a\ contractible\ circle\ }C_{01}\mathrm{\ but\ }r_{10}\mathrm{\ does\ not}\\
0 & \mathrm{\ neither\ }r_{10}\mathrm{\ nor\ }r_{01}\mathrm{\ has\ a\ contractible\ circle}\\
\end{array}\right.
\end{equation}
If $r'$ is a resolution with $|r'| = |r| + 2$ but no pair of crossings $c_{1}$ and $c_{2}$ exists with $r' = r_{11}$ as above then $\pairing{r}{r'} = 0$.
\end{defn}

\noindent This yields an explicit description for the chain complex supported on noncontractible circles. for each $b \in \glyph{\Sigma}{B}$ let
$$
\lk{O}_{n,i}(T,b) = \big\{\,r\, \big|\, r \in  \lk{O}_n(T,b),\,|r|=i,\,\beta(r)=b\,\big\}
$$
The deformed chain complex consists of the vector spaces spanned by $\lk{O}_{n,i}(T,b)$
$$
CT_{i}(T,b):= Span_{\field{T}} \big \{ r \in \lk O_{n,i} (T,b)  \big \} 
$$
equipped with the boundary map 
$$
\partial_{n,i}: CT_{i}(T,b) \longrightarrow CT_{i+2}(T,b) 
$$
defined by
$$
\partial_{n,i}(r)= \sum_{r' \in  \lk O_{n,i+2}(T,b)} \pairing{r}{r'}\,r'
$$  
Since this complex is homotopy equivalent to $\widetilde{SK}(T,b)$, theorem \ref{theorem:main} implies

\begin{theorem}
Let $T$ be diagram of $I$-bundle tangle $\lk T$ in $\Sigma \times (-1,1)$ with boundary $B \subset \partial \Sigma$. For each glyph $b \in \glyph{\Sigma}{B}$, let $CT_{*}(T,b)= \bigoplus_{i \in Z} CT_{i}(T,b)$ and  let $\partial_{n}$ be the boundary map defined above. Then the (stable) chain homotopy class of $(CT_{*}(T,b), \partial_{n})[-n_{+}(T)]$ 
is preserved by the Reidemeister I,II, and III moves and thus is an invariant of $\lk{T}$. 
\end{theorem}

\begin{defn}
The $i^{th}$ homology group of  $(CT_{*}(T,b), \partial_{n})$ will be denoted $HT_{i}(T,b)$. The homology groups of 
$(CT_{*}(T,b), \partial_{n})[-n_{+}(T)]$ will be denoted $\underline{HT}_{*}(\mathcal{T},b)$. Finally let 
$$\underline{HT}_{*}(\mathcal{T})  = \displaystyle{\bigoplus_{b \in \glyph{\Sigma}{B}} \underline{HT}_{*}(\mathcal{T},b)}$$ 
\end{defn}

\section{Behavior under symmetry}

\noindent If $T$ is a tangle diagram for $(\Sigma, B)$, the {\em mirror} of $T$ is the diagram where all the right handed crossings are changed to left-handed crossings and vice-verse. Alternately, it is the tangle $T$ considered in $(\overline{\Sigma}, B)$. $\overline{T}$ is also the diagram for the image of $\mathcal{T}$ under the orientation preserving map $\Sigma \times (-1,1) \rightarrow \overline{\Sigma} \times (-1,1)$ induced by $(x,t) \rightarrow (x, -t)$. \\
\ \\
\noindent In addition, if $g \in \glyph{\Sigma}{B}$ then $-g$ is the glyph with the same configuration of arcs and circles, but such that the label on a circle $C$ is the negative of that which $g$ assigns to $C$.

\begin{prop}
For each $g \in \glyph{\Sigma}{B}$ there is an isomorphism $\underline{HT}_{\delta}(T,-g) \cong \underline{HT}_{\delta - 2k(g)}(T,g)$
\end{prop}

\noindent{\bf Proof:} Suppose $r$ is a state contributing to $\underline{CT}(T,-g)$. Let $\Phi(r) \in \underline{CT}(T,g)$ be the state found by changing the signs on all the circles, but preserving the choice of resolutions at each crossing. $\partial\, r$ is built out of two operations: 1) cleaving a contractible circle from a noncontractible circle and reattaching it (possibly to  itself to get an annulus), or 2) joining the boundary circles of an annulus that are labeled with a $+$ and a $-$ to get a contractible circle, and then reattaching the contractible circle (possibly to itself). Running through the possible sign assignments shows that changing the signs on all the circles preserved the boundary. Thus, $\Phi$ is a chain map. Since $\Phi$ is also an isomorphism, we have established the result up to the change in gradings. The $\delta$-grading of $r$ is $2h(r) - q(r) - n_{+}(T) = 2|r| - (|r| + k(g)) - n_{+}(T) = |r| - k(g) - n_{+}(T)$. For $\Phi(r)$ the $\delta$-grading is $|r| - k(-g) - n_{+}(T) = |r| + k(g) - n_{+}(T) = \delta(r) + 2k(g)$. The statement in the proposition follows.  $\Diamond$

\begin{prop}
Let $\overline{T}$ be the mirror of $T$ in $(\Sigma, B)$. For each $g \in \glyph{\Sigma}{B}$ there is an isomorphism $\underline{HT}_{\delta}(\overline{L},g) \cong \underline{HT}^{-\delta}(L,-g)$. 
\end{prop}

\noindent{\bf Proof:} Let $r$ be a state for $T$ in grading $\delta$, and define $\Phi(r)$ to be the state for $\overline{T}$ which results in the same resolution diagram, but assigns the opposite signs to circles. Thus $\Phi(r)$ resolves the crossing $c$ using the opposite code from $r$: $\Phi(r)(c) = 1 - r(c)$ and if $\gly{r} = g$ then $\gly{\Phi(r)} = -g$. Finally, let $f_{r}$ be the element in the dual vector space $\underline{CT}^{\delta}(L,g)$ which evaluates to $1$ on $r$ and $0$ on all other generators. Then if $r'$ is a state for $T$ is $\delta$-grading $\delta-2$, we have
$\delta(f_{r})(r') = f_{r}(\partial(r')) = \sum_{r'} \pairing{r'}{r}$. Thus,
$$ \delta(f_{r}) = \sum_{r'} \pairing{r'}{r}\,f_{r'} $$
where the sum is taken over those states found by changing two $1$-resolved crossings in $r(T)$ to $0$-resolved crossings. On the other hand, due to the change in resolution codes, and the indifference to changing signs on all the circles, the same formula describes $\partial_{\overline{T}}(\Phi(r))$. Thus $\Phi$ is a chain isomorphism, and we turn to deciphering the effect on $\delta$-gradings. The grading $\delta$ of $f_{r}$ is $2h(r) - q(r) - n_{+}(T) = h(r) - k(g) - n_{+}(T)$. The grading for $\Phi(r) = (n-h(r)) - k(-g) - n_{+}(\overline{T})$ $= (n-h(r)) + k(g) - n_{-}(T)$ $= (n-n_{-}(T)) - h(r) + k(g) = n_{+}(T) - h(r) + k(g)$, but this is just $-\delta$, and the statement in the proposition follows. $\Diamond$\\
\ \\
\noindent Furthermore, since the coefficients are take in the field $\field{T}$, the universal coefficient theorem ensures that $\underline{HT}^{-\delta}(T,-g) \cong \underline{HT}_{-\delta}(T,-g)$, and thus, using both propositions, 
$$\underline{HT}_{\delta}(\overline{T},g) \cong \underline{HT}_{-\delta}(T,-g) \cong \underline{HT}_{-\delta - 2k(g)}(T,g)$$

\section{Thomas C. Jaeger's Trick}\label{sec:Jaeger}

\diagram[0.7]{jaeger}{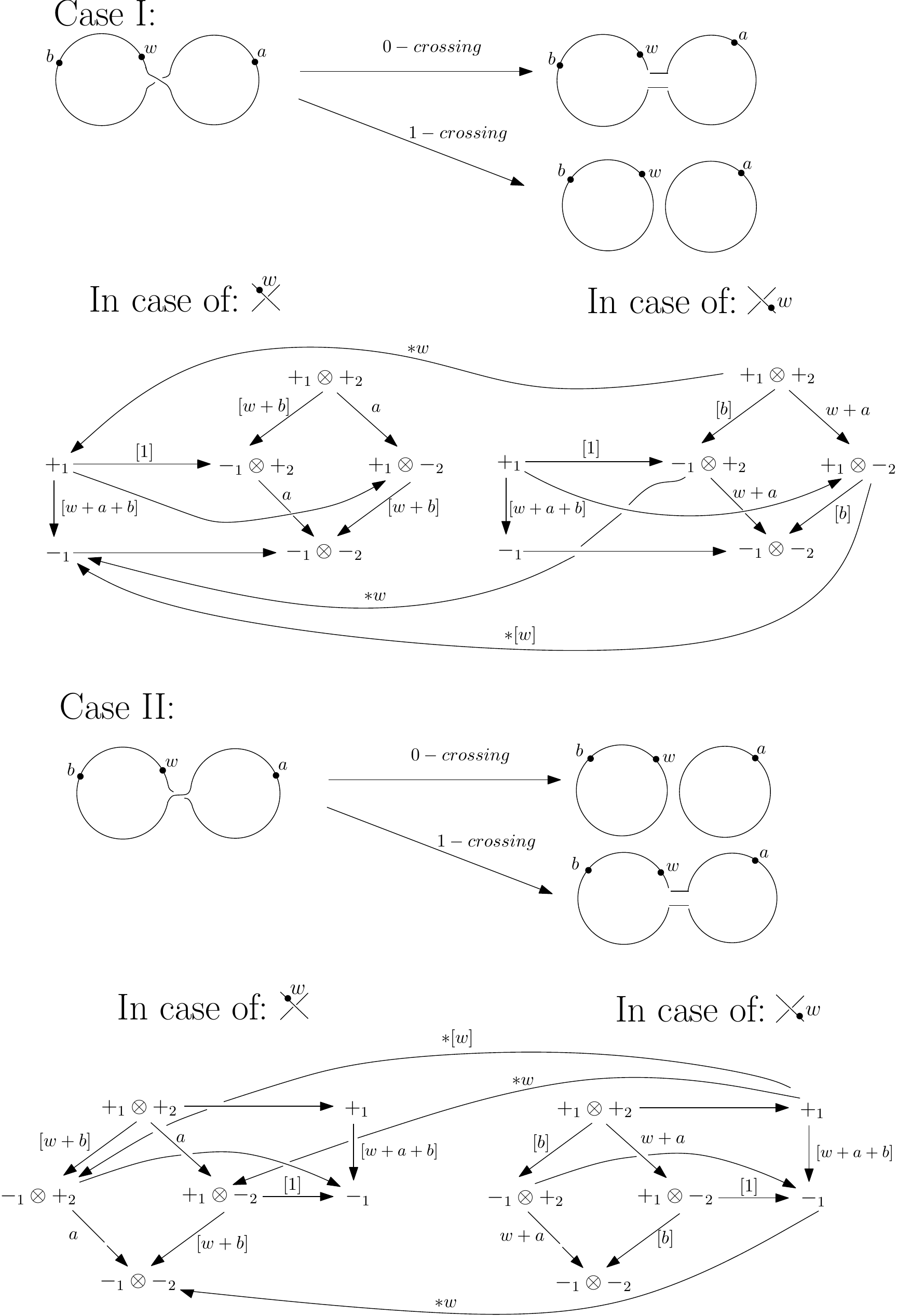}{Illustrations of T. Jaeger's chain isomorphisms. The maps $\cdot w$ are added to the horizontal identification map $I$ to create the chain isomorphism $\Psi$ going from left to right. By $[a]$ we mean either the formal variable $a$ if the corresponding object is a contractible circle, or $0$ if the object is not a contractible circle . Similarly, the Khovanov map $[1]$ will be the normal Khovanov map if the left object is a circle appearing in the $0$-resolution of Case $I$ or the $1$-resolution of Case $II$. Otherwise, the map is zero when the left object is an arc or noncontractible circle.} 

\noindent Before proving Reidemeister invariance, we follow \cite{TJ} and prove a lemma that allows us to move the weight on an edge past a crossing. This will greatly simplify the proof of invariance in section \ref{sec:invariance}. The proof consists of adapting that in \cite{TJ} to arcs and noncontractible circles. The argument will be presented for a weight moving on the over-strand of a crossing as in \ref{fig:jaeger}. \\
\ \\
\noindent Let $c$ be the crossing the weight $w$ moves past, $L$ be the diagram with edge weights before the weight is moved, and $L'$ be the diagram after the weight is moved. Let $\overline{D}_{c}$ be the map taking a resolution $r$ for $L$ to 1) $0$ if $r(c) = 0$ and, 2) the state(s) for $L'$ where $r(c') = r(c')$ for $c \neq c'$, $r(c) = 1$, and the signs on the circle(s) abutting $c$ computed using the Khovanov maps $\partial_{Kh}$. Let
$I$ be the map which takes the state $r$ for $L$ to the state $r$ for $L'$. 

\begin{prop}
The map $\Psi = I + w\cdot \overline{D}_{c}$ is a $\delta$-graded, chain isomorphism from $\underline{SK}(L)$ to $\underline{SK}(L')$. 
\end{prop}

\noindent{\bf Proof:} 
First we verify that $\overline{D}_{c}$ preserves the $\delta$-grading. If $r$ occurs in bi-grading $(h(r),q(r))$ then the sum of the decorations for $r$ is $q(r) - h(r)$. The homological grading of $\overline{D}_{c}(r)$ is $h(r) - 1$ since we change a $1$ resolution to a $0$ resolution. The sum of the signs decreases by $1$ since $\overline{D}_{c}$ comes from the Khovanov maps. Thus $\overline{D}_{c}(r)$ occurs in  bi-grading $(h(r) - 1, (q(r) - h(r) - 1) + (h(r) - 1)) = (h(r)-1, q(r) - 2)$. Thus, $\delta(\overline{D}_{c}(r) = 2(h(r) - 1) - (q(r) - 2) = 2h(r) - q(r) = \delta(r)$. \\
\ \\
The boundary map for $\underline{SK}(L)$ can be written as a sum $\partial_L = D_{L}$ + $\partial^{Kh}_{L,c}$ + $E_{L}$ + $E_{L, c}$ where $D_{L}$ is the sum of the Khovanov maps coming from changing the resolution from $0$ to $1$ at $c' \neq c$, $\partial^{Kh}_{L,c}$ is the Khovanov map for the resolution change at $c$, $E_{L}$ is the sum of the vertical maps for contractible circles not abutting $c$, and $E_{L,c}$ is the sum of the vertical maps for the circles abutting $c$. We can do this for both $\underline{SK}(L)$ and $\underline{SK}(L')$. As the effect of $\overline{D}_{c}$ is that of a Khovanov map (for the mirror crossing at $c$) we know that $D_{L'}\overline{D}_{c} = \overline{D}_{c}D_{L}$ and $E_{L}\overline{D}_{c} = \overline{D}_{c}E_{L}$. Likewise, $D_{L} = D_{L'}$ and $E_{L} = E_{L'}$ as nothing has changed in the diagrams outside of a small region of $c$. Furthermore, $\partial^{Kh}_{L,c} = \partial^{Kh}_{L',c}$ since these do not depend on the weights. To see that
$\big(I + w\cdot \overline{D}_{c}\big)\partial_{L} = \partial_{L'}(I + w\cdot \overline{D}_{c}\big)$ we need to see that
\begin{equation}\label{eqn:ident}
E_{L,c} + w\cdot \overline{D}_{c}\partial^{Kh}_{L,c} + w\cdot\overline{D}_{c}E_{L,c} = E_{L',c} + w\cdot \partial^{Kh}_{L',c}\overline{D}_{c} + w\cdot E_{L',c}\overline{D}_{c}
\end{equation}
\ \\
\noindent First we show that $\overline{D}_{c}E_{L,c} = E_{L',c}\overline{D}_{c}$. Note that both will be zero unless $r(c) = 1$, as $E_{L,c}$ does not change the resolution of $c$. Furthermore, $E_{L,c}$ is zero unless there is a contractible circle decorated with a $+$ in the image of $\overline{D}_{c}$. Consequently, the right term is $0$ when the object(s) abutting $c$ are 1) two arcs, 2) two noncontractible circles (merging them will be $0$ or a $-$ decorated circle), 3) a circle and an arc, 4) a circle and a noncontractible circle, 5) a single arc, or 6) a single noncontractible circle. $E_{L,c}$ will be zero on 1),2),5),6). For 3) and 4) $E_{L,c}$ may not be zero, but it will leave a $-$ decorated circle, and the merging of a $-$ circle into an arc or a noncontractible circle will be $0$. The remaining cases are those where $c$ only abuts one or two contractible circles. If there are two contractible circles, both sides will be zero unless both circles are decorated with $+$'s. In that case
the image of $E_{L',c}\overline{D}_{c}$ will be $(a + w + a')\cdot (-)$ where $-$ adorns the circle resulting from merging the two circles at $c$. $E_{L,c}$ will have image $(a+w)(-\otimes + ) + a'(+ \otimes -)$, but $\overline{D}_{c}$ will map both generators to $(-)$ giving the desired equality.
If there is only a single circle, it must have a $+$ decoration for both sides to be non-zero. In this case, the right side will yield $E_{L',c}\big((+\otimes -) + (-\otimes +)\big) = a(- \otimes -) + (w + a')(- \otimes -) = (a + w + a')(- \otimes -)$ while the left side gives
$\overline{D}_{c}\big((a+w+a')(-)\big) = (a + w + a')(- \otimes -)$ since all the weights will be on the single circle. \\
\ \\
\noindent Canceling these terms from \ref{eqn:ident} leaves us with the task of verifying  
$$
E_{L',c} + E_{L,c} =  w\cdot \left(\overline{D}_{c}\partial^{Kh}_{L,c} + \partial^{Kh}_{L',c}\overline{D}_{c}\right)
$$
We first describe the left side of this equation. It is zero unless a circle abutting $c$ is contractible and decorated with a $+$.  If there are two objects abutting $c$, the image will be the sum over the zero, one, or two contractible $+$ circles where we change the $+$ to a $-$ and multiply by $w$
(since the each arc has $w$ in one and only one of $L'$ or $L$, and the other weights don't change). On the other hand if there is one object
abutting $c$, the image of the left side will be zero except when this is a contractible $+$-circle, and then both terms map it to $(a+w+a')(-)$. The left side is thus $0$ always in this case. \\
\ \\
\noindent The right side can likewise be computed directly. Only one of the compositions $\partial^{Kh}_{L',c}\overline{D}_{c}$  or $\overline{D}_{c}\partial^{Kh}_{L,c}$ will be non-zero depending upon $r(c)$. These compositions start at $r$ (for $L$) and map to $r(L')$ with different decorations. If there are two objects abutting $c$, the image will be the sum over the $+$ contractible circles, changing the decoration on that circle to $-$. When we multiply by $w$ we get the same image as $E_{L',c} + E_{L,c}$. If there is only one object, the image is $0$ unless that object is a contractible $+$ circle, then the image is $(+) \rightarrow (+-) + (-+) \rightarrow (--) + (--) = 0$ oe $(-) \rightarrow (--) \rightarrow 0$. In either case, the two sides are equal. Thus $\Psi$ is a chain map. \\
\ \\
\noindent To see that it is a chain isomorphism, consider $\Xi = I^{-1} + w\cdot \underline{D}_{c}$ where $\underline{D}_{c}$ is defined identically to $\overline{D}_{c}$ but mapping from $L'$ to $L$. Then $\Xi\Psi = \mathbb{I} + w\cdot \big(\underline{D}_{c}I + I^{-1}\overline{D}_{c}\big) + w^{2}\big(\underline{D}_{c}\overline{D}_{c}\big)$. However, the image of $\overline{D}_{c}$ is supported on those resolutions with $r(c) = 0$, while the domain of $\underline{D}_{c}$ is limited to those with $r(c) = 1$. Consequently, the last term vanishes. Furthermore, $\underline{D}_{c}$ and
$\overline{D}_{c}$ have the same affect on a decorated state. All that differs is whether the image and domain are taken to be in the complex for $L$ or that for $L'$. It is straightforward to verify that $I\underline{D}_{c}I = \overline{D}_{c}$. Thus $\Xi\Psi = \mathbb{I}$, and reversing the roles of $L'$ and $L$ shows that $\Psi\Xi = \mathbb{I}$ as well. Thus $\Psi$ is a chain isomorphism. $\Diamond$.\\
\ \\
\noindent{\bf Note:} The argument above does not presume that $c$ is a certain type of crossing. However, it does require that
the weight $w$ before its transit occurs on one circle, and after its transit occurs on the other circle, in the case of two contractible circles abutting the crossing $c$. Since the circles can occur to the left and right, or to the top and bottom, of the crossing $c$ the weight must move along its link component. It cannot jump from one component to the other. However, both of the following transits are allowed: 
$$
\inlinediag{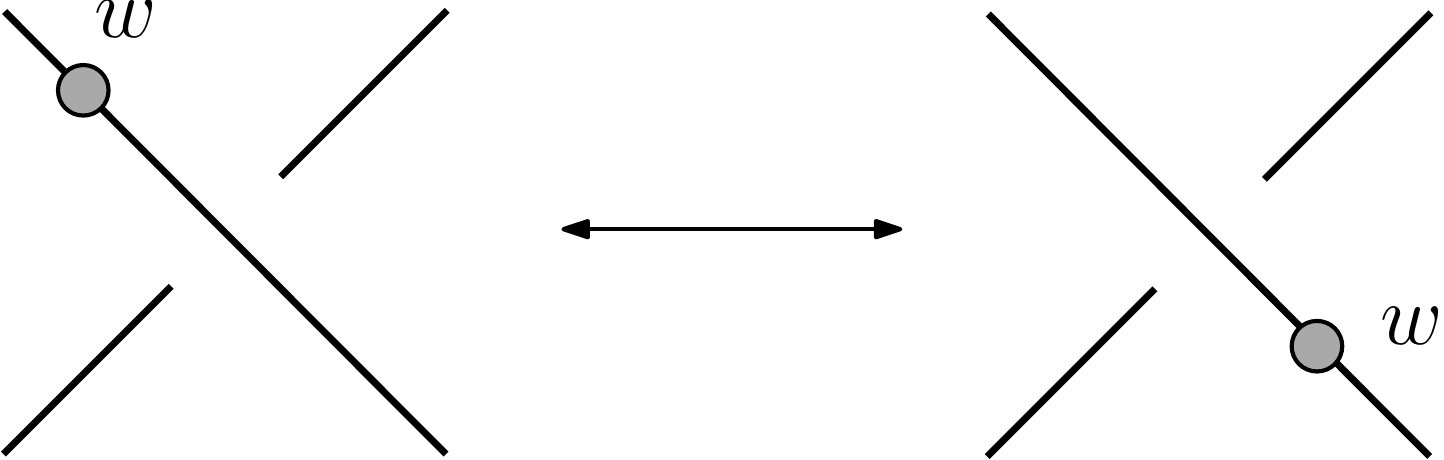} \hspace{0.5in} \inlinediag{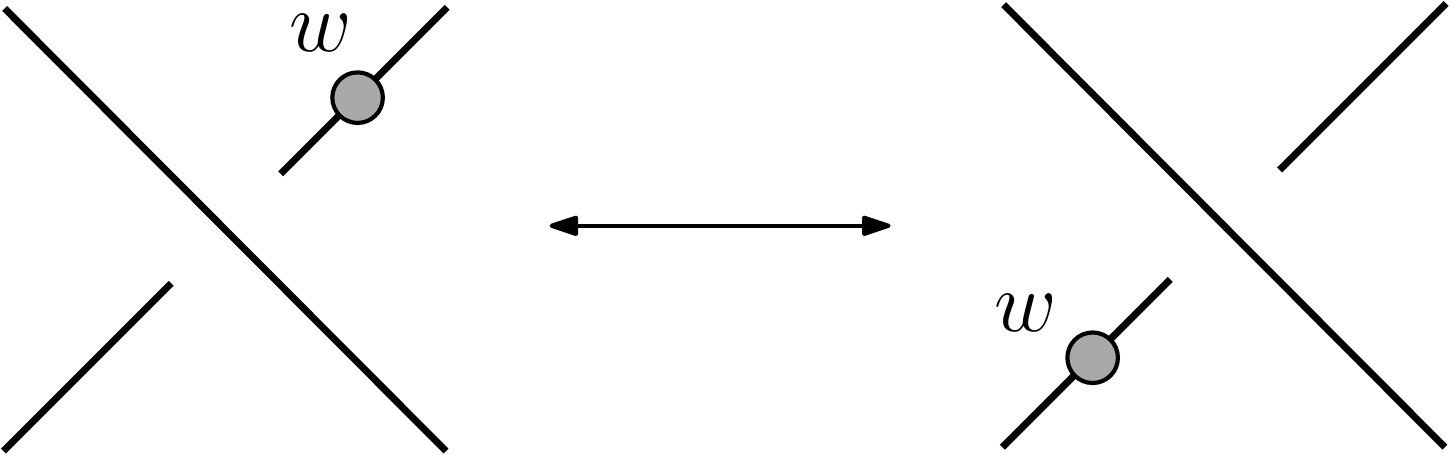} 
$$
\ \\
\noindent In some cases the twisted skein homology is an effective substitute for the corresponding Khovanov skein homology. Following T. Jaeger's argument in \cite{TJ} we can use this weight shifting isomorphism to prove
\begin{theorem}
If $T$ is a tangle consisting only of arcs in $(\Sigma,B)$, then $\underline{SK}(T)$ is isomorphic to the the Khovanov skein homology defined in \cite{APS2}.
\end{theorem}
\noindent We simply move all the weights along their components in $T$ until they are at a boundary point in $B$. Note that this theorem applies to reduced Skein homologies when $T$ is equipped with a marked point on each component. One removes a neighborhood of the marked point to obtain a tangle in $\Sigma'$ without closed components whose Khovanov skein homology is isomorphic to the reduced skein homology for $T$ in $\Sigma$  determined by those marked points. The isomorphisms respect the decompositions of the respective homologies along glyphs. Note that we have used a slightly different convention for the quantum gradings, and for the Khovanov maps, and the theorem presumes that one modifies either the construction in \cite{APS1} or ours accordingly. \\
\ \\
\noindent{\bf Example: cf. \cite{TJ}:} If $T$ is a tangle in $(D^{2},P_{2n})$, where $P_{2n}$ is a set of $2n$ points in $\partial D^{2}$, then the twisted homology is isomorphic to the skein homology. Glyphs here are crossingless matchings of the points in $P_{2n}$. 

\section{Invariance of the twisted homology under Reidemeister moves}\label{sec:invariance}

\noindent To prove $\underline{SK}(T)$ is invariant under the Reidemeister moves, we will first use Thomas Jaeger's isomorphism to move the weights to the bottom of the local diagram, see Figures \ref{fig:RI} and \ref{fig:RIIweights}. Once we have done that, the result follows from a modification of the standard proof for the invariance of Khovanov homology. We give an account of the modifications for each of the three types of moves below. As in \cite{LR}, there is an important technical point to note: the chain complexes before and after a Reidemeister move will be defined over different fields. The details for relating the complexes through stable isomorphism are in \cite{LR}, but as this causes no serious trouble they will be omitted here.

\diagram[0.7]{RI}{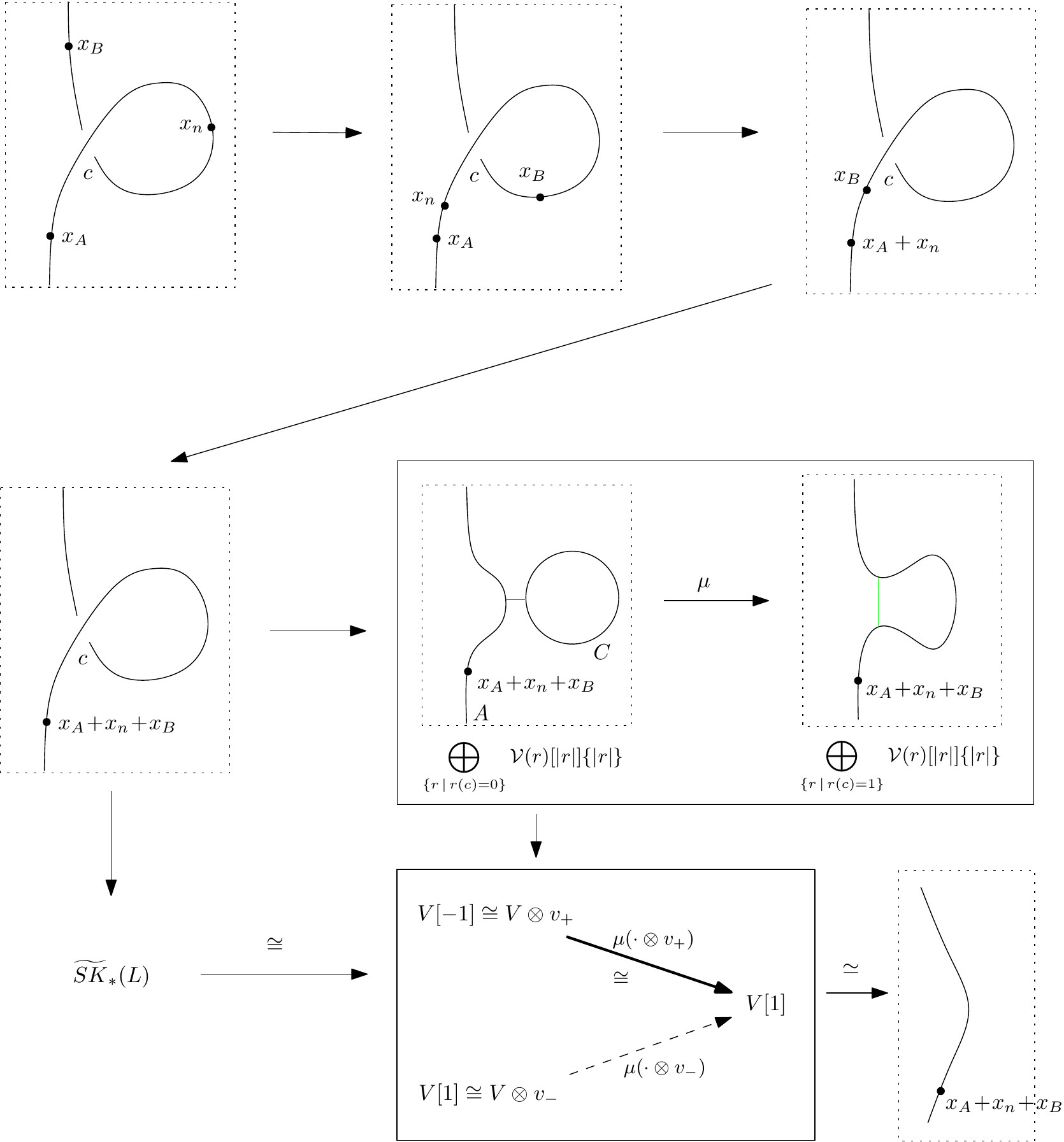}{A schematic representation of the proof of invariance under the first Reidemeister move, for a positive crossing. In the top row, we use Thomas Jaeger's weight shift isomorphisms to move all the weights to the bottom of the diagram. Resolving the crossing $c$ in both ways gives a finer view into the complex. Regardless of whether the local arc is on a global arc, a noncontractible circle, or a noncontractible circle, the thickened arrow is an isomorphism (however, the nature of the dashed arrow depends upon the type of the local arc. When the complex is reduced along the thickened arrow we obtain the complex for the diagram after the RI move, but with the weight $x_{A} + x_{n} + x_{B}$ on the local arc.}

\noindent{\bf Invariance under Reidemeister I moves:} Figure \ref{fig:RI} shows the complex for a diagram prior to and after an RI move applied to a right-handed crossing.  
Let $T$ be the diagram before the RI move, and let $T'$ be diagram afterwards. For a sake of simplicity, the $\delta$-graded chain complex $\widetilde {SK}(L')$ will be called $V$. As described in \cite{LR}, we decompose $\widetilde{SK}(T)$ into $(\oplus_{{r(c)=0}} \lk V(r)) \oplus (\oplus_{r(c)=1} \lk V(r))[1]$. When $r(c) = 0$, there is a local contractible circle $C$ in the diagram. Using the decoration on $C$ we further decompose $\oplus_{r(c)=0} \lk V(r) =( V \otimes v_+) \oplus (V \otimes v_-) \cong V[-1] \oplus V[1]$ where $[\cdot]$ is the shift in $\delta$-grading ($\otimes v_{+}$ increase the quantum grading by $1$, and thus {\em decreases} the $\delta$-grading by $1$). since we have move the weights, the complex corresponding to this circle is $0 \rightarrow \field{T} \stackrel{0}{\rightarrow} \field{T} \rightarrow 0$, just as it is for Khovanov homology.
On the other hand $\oplus_{r(c)=1} \lk V(r) \cong V[1]$.\\
\ \\
Changing the resolution at $c$ from $0$ to $1$ corresponds to a term in the boundary map. Since this resolution change merges the a $+$ contractible circle the map $\mu(.\otimes v_{+}):$ $V \otimes v_+ \longrightarrow \bigoplus_{c=1} \lk V(L_S)$ is an isomorphism. This is true regardless of whether the local arc belongs to a global arc, a contractible circle, or a noncontractible circle. As this is known for contractible circles, it is enough to check if it is true when our local arc is  part of a global arc or a noncontractible circle. Looking back at Figures \ref{fig:Khovanov-map2} and \ref{fig:Khovanov-map3}, we see that this is indeed the case.\\
\ \\
\noindent We can use Gaussian elimination with this isomorphism to cancel $V \otimes v_+ \cong V[-1]$ and $\bigoplus_{c=1} \lk V(L_S) \cong V[+1]$. The result is a chain homotopy equivalent chain complex supported on $V \otimes  v_- \cong V[1]$ (see \cite{LR}, section 5.2). Since the vertical differential from $C$ is trivial, there {\em is no alteration} to the boundary map on $V \otimes  v_{-}$. Indeed,
if the local arc is part of a global arc or of a noncontractible circle, then $\partial_{new}=\partial|_{(V \otimes v_{-})}$ since $\mu (v_{+,1} \otimes v_{-,c})=0$ and $\mu (v_{0} \otimes v_{-,c})=0$. Thus the cancellation introduces no new terms in the boundary map.\\
\ \\
\noindent If the arc in local diagram is a part of a contractible circle $D$, we have $\partial_{new} (v \otimes v_{+} \otimes v_{-,c})= \partial|_{(V \otimes v_{-})} \big(v \otimes v_{+} \otimes v_{-,c}+[C].(v \otimes v_{-}\otimes v_{-,c})\big)$, where the second and third factors stand for decoration on $D$ and $C$ respectively. But since there is no formal variable on $C$, $[C]=0$. As a result, $\partial_{new}=\partial$ on $V \otimes v_-$. \\
\ \\
\noindent There is a difference between $((V \otimes v_-), \partial_{new})$ and  $(\widetilde{SK}(T')[1], \partial')$, however. The local arc in $V \otimes v_{-}$ is decorated with $x_{A} + x_{n} + x_{B}$ while that for $L'$ will be decorated with some formal variable $y_{i}$. The map $y_{j} \rightarrow x_{j}$, $y_{i} \rightarrow x_{A} + x_{n} + x_{B}$ will induce an inclusion $\field{T'} \hookrightarrow \field{T}$. Using this inclusion, $(\widetilde{SK}(T')[1], \partial') \otimes \field{T}$ is chain isomorphic to $((V \otimes v_-), \partial_{new})$. Furthermore, since $n_{+}(L)=n_{+}(L')+1$ 
$$
((V \otimes v_-) [-n_{+}(T)], \partial_{new}) \cong (\widetilde {SK}(T')[1][-n_{+}(T')-1], \partial') \cong (\widetilde {SK}(T')[-n_{+}(T')], \partial') $$  
So $(\widetilde{SK}(T)[-n_{+}(T)],\partial) \cong (\widetilde {SK}(T')[-n_{+}(T')], \partial')$. Thus an RI move on a right-handed crossing in $T$ induces a stable chain homotopy equivalence $\underline{SK}(T) \simeq \underline{SK}(T')$ of $\delta$-graded chain complexes.  \\
\ \\
\noindent The same approach works when a Reidemeister I move is applied to a left-handed crossing $c$. Once we move the weights past the crossing we obtain a complex where the $1$ resolution at $C$ has a local contractible circle $C$. The change in resolution introduces a division of the arc, and it suffices to show that $\Delta: V \longrightarrow V_+[1] \bigoplus V_-[1]$ followed by projection onto $V_{-}[1]$ is surjective, where $\bigoplus_{r(c) = 1} \lk V(L_S) \cong V_{+}[1] \bigoplus V_{-}[1] \cong V \bigoplus V[2]$. However, dividing a circle from a contractible circle, arc or noncontractible circle always has a term $v \rightarrow v \otimes v_{C,-}$ (see Figures \ref{fig:Khovanov-map1}, \ref{fig:Khovanov-map2}, and \ref{fig:Khovanov-map3}). So the same cancellation procedure applies, and a similar argument shows that $\widetilde{SK}(L)$ is chain homotopy equivalent to $V \cong V_{+}[1]$. Once again, since there are no weights on $C$, the new boundary map on $V$ is identical to the restriction of the boundary map from $\widetilde{SK}(L)$. The gradings are the same as those for $\widetilde{SK}(L')$, but the weights on the local arcs are again $x_{A} + x_{n} + x_{B}$ and $y_{i}$. The inclusion induced by $y_{i} \rightarrow x_{A} + x_{n} + x_{B}$ will provide the stable isomorphism as before. 

\diagram[0.75]{RIIweights}{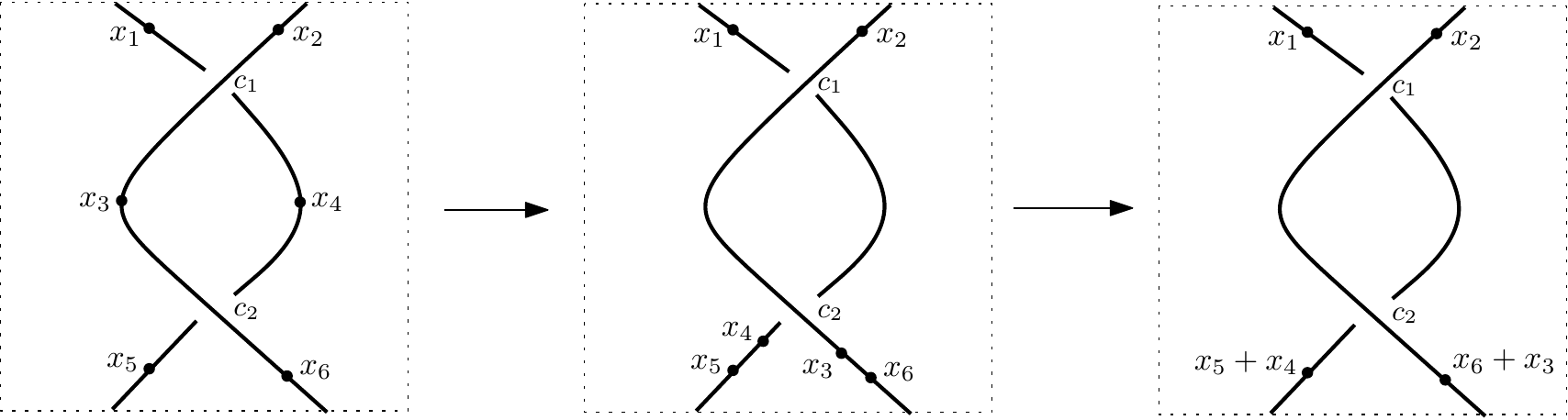}{}
\diagram[0.75]{RIIShift}{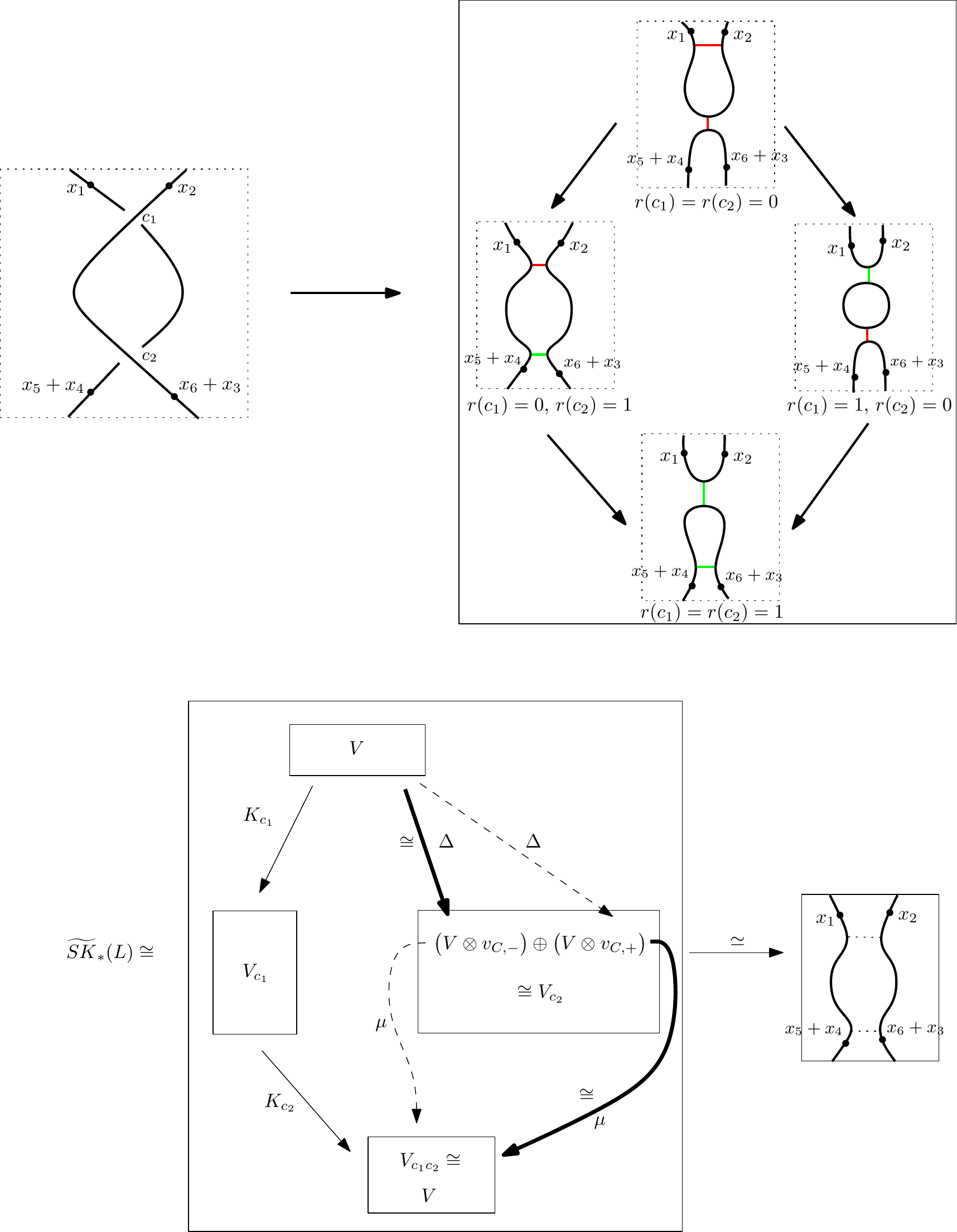}{A schematic representation of the proof of invariance under the second Reidemeister move. The complex is again reduced along the thickened arrows, which are both isomorphisms regardless of whether the local arcs lie on global arcs, contractible circles, or noncontractible circles. The dashed arrows do depend upon the type of the local arcs, however. Nevertheless, if we cancel the bottom arrow first, and then the top arrow, we introduce no new terms in the boundary map since $[C]=0$. These cancellations produce the chain homotopy equivalence complex in the lower right.}  

\noindent{\bf Invariance under Reidemeister II moves:} We start with a diagram as in the left side of Figure \ref{fig:RIIweights}. Once again we shift all the weights to the bottom of the diagram. Notices that $x_{1}$ and $x_{2}$ will be present throughout, and can be ignored until the very end. As before $T$ will be the tangle before the move and $T'$ will be diagram after using an RII move to remove crossings $c_1$ and $c_2$. We can then decompose the complex into parts based on the resolutions at $c_{1}$ and $c_{2}$. Doing this gives the diagrams and complexes depicted in Figure \ref{fig:RIIShift}. Using the cancellations described by Figure \ref{fig:RIIShift} In produces a homotopy equivalent chain complex which can be described as $\lk V_{c_1}$ with a perturbed boundary map $\partial_{new}$. As described in \cite{LR}, $\partial_{new}=\partial_{\lk V_{c_1}}+ [C]K_{c_1}K_{c_2}$ where $[C]$ is the sum of formal variables on $C$. since we moved the weights $[C] = 0$. Thus $(\widetilde{SK}(L),\partial) \cong (\lk{V}_{c_1}, \partial_{\lk{V}_{c_1}})$. \\
\ \\
\noindent $(\lk{V}_{c_1}, \partial_{\lk{V}_{c_1}})$ is almost identical with $\widetilde{SK}(T')[1]$ (where the shift comes from the resolution at $c_{1}$), except that the weights on the two local arcs in $T'$  will be $y_{L}$ for the left arc and $y_{R}$ for the right arc.  In $\lk{V}_{c_1}$ the corresponding weights are  
$x_{1} + x_{4} + x_{5}$ and $x_{2}+ x_{3} + x_{6}$. As before the map $y_{L} \rightarrow x_{1} + x_{4} + x_{5}$, $y_{R} \rightarrow x_{2}+ x_{3} + x_{6}$, $y_{j} \rightarrow x_{j}$ can be used to find a chain isomorphism between $\widetilde{SK}(T')[1] \otimes \field{T}$ and $(\lk{V}_{c_1}, \partial_{\lk{V}_{c_1}})$. Finally, $n_{+}(T)=n_{+}(T')+1$ so $$(\widetilde{SK}_{\ast}(T), \partial)[-n_{+}(T)] \simeq (\lk{V}_{c_1}, \partial_{\lk{V}_{c_1}})[-n_{+}(T)] \simeq (\widetilde {SK}_{*}(T'), \partial')[1][-n_{+}(T')-1]$$
or $(\underline{SK}_{\ast}(T), \partial) \simeq (\underline{SK}_{*}(T'), \partial')$

\diagram[0.7]{RIII_LShift}{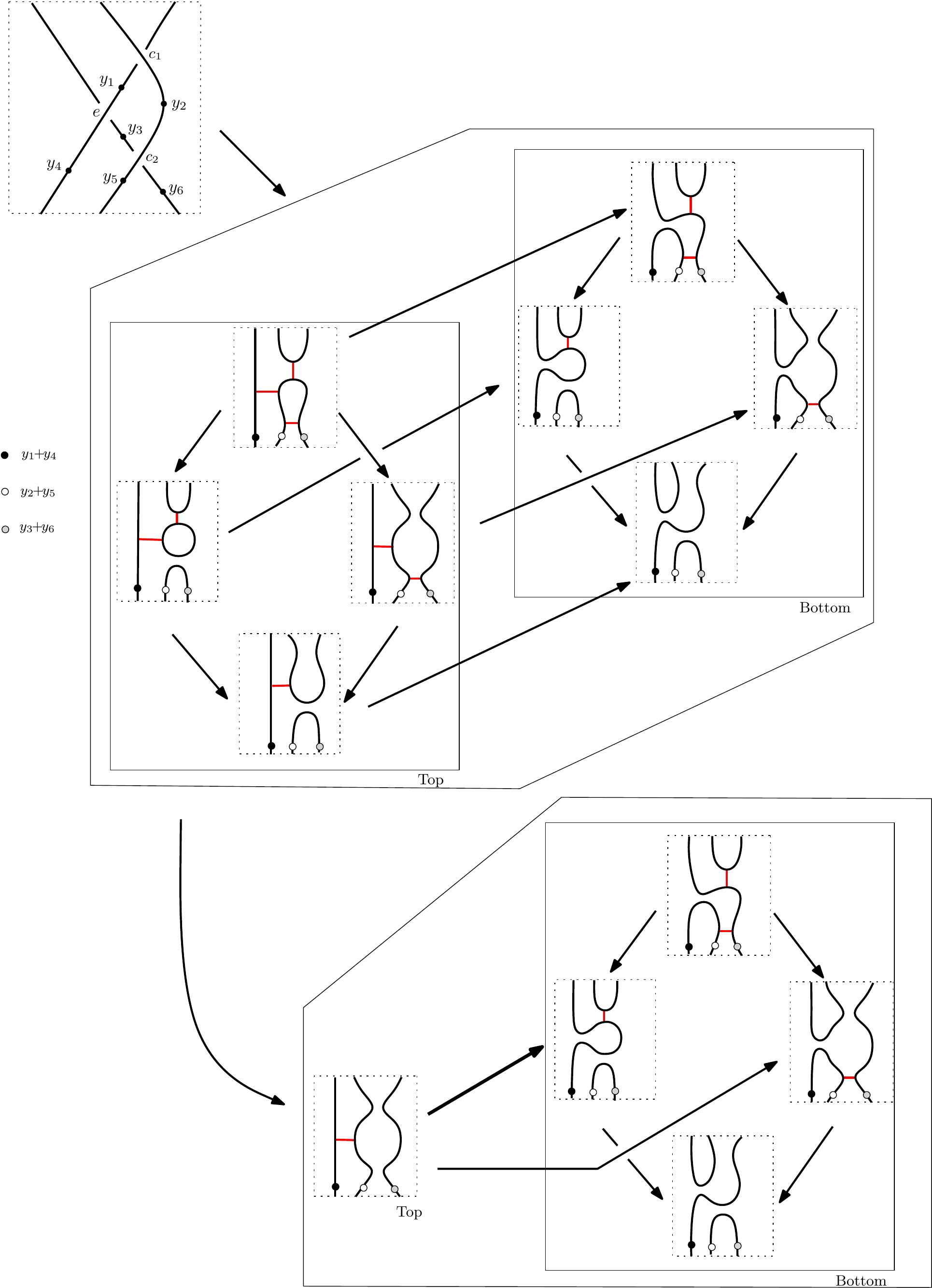}{The local picture for a diagram before an RIII move. We decompose the complex along the
eight possible ways of resolving the local crossings. The four resolutions with crossing $c$ resolved using a $0$ resolution
replicate the diagrams in the proof of RII invariance. Canceling along the same isomorphisms, in the top diagram, gives the
diagram at the bottom. A new term in boundary map may occur from the thicker arrow in the bottom picture; however, it replicates the map for changing the resolution at $c_{2}$. The colored circles represent the weights, listed on the left.}
\diagram[0.7]{RIII_RShift}{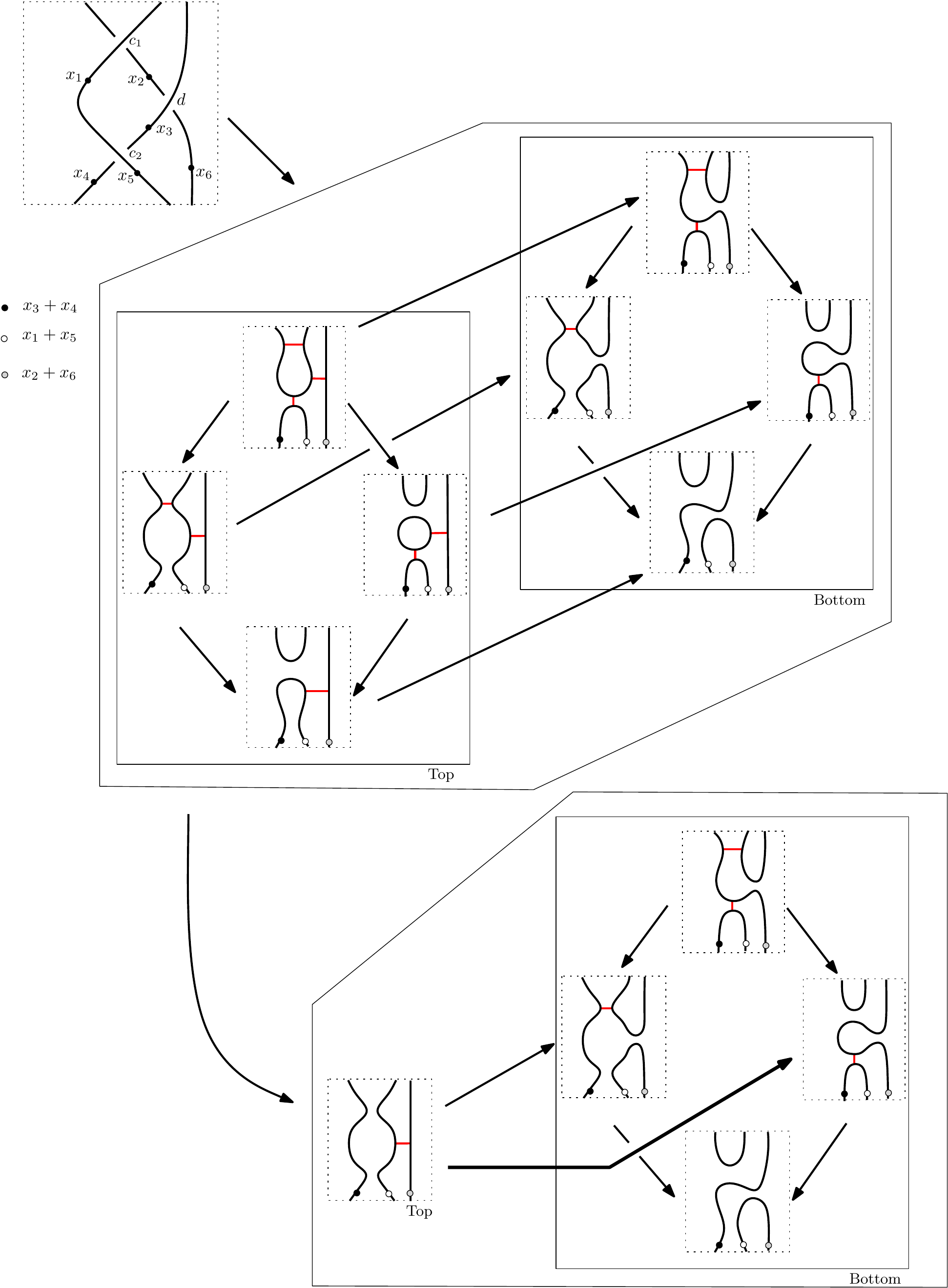}{The local picture for a diagram after the RIII move. Once again there is a new term in boundary map,shown by the thicker arrow in the bottom picture. This time it replicates the map for division at $c_{1}$.}

\noindent{\bf Invariance under Reidemeister III moves:} The strategy is the same as for RI and RII moves:
\begin{enumerate}
\item Move all the weights to the bottom of the local picture
\item Decompose the chain complex according to the resolutions of the crossings in the local diagrams
\item Locate {\em contractible} circles in the diagram with $0$ formal weight, and further decompose based on the decoration of these circles
\item Use the merging of $+$ contractible circles, and the division of of an object into the same type object and a $-$-contractible circle to identify isomorphism on the summands in the decomposition
\item by canceling these isomorphisms as in the untwisted Khovanov homology, show that one can obtain a chain homotopic complex where the crossings have been simplified
\end{enumerate}
We depict this process for two local diagrams related by an RIII-move in Figures \ref{fig:RIII_LShift} and \ref{fig:RIII_RShift}. We assume the weights for the edges not shown in the picture are $z_{7}, \ldots z_{n}$. The bottom diagrams in the last step of Figures \ref{fig:RIII_LShift} and \ref{fig:RIII_RShift} are the same, up to the weights, and the maps from the top diagram to the lower diagram are purely Khovanov maps (since the weight on the contractible circles is $0$) and are known to be the same, \cite{Bar2}. Call these complexes, with the weights depicted, $R_{1}$ and $R_{2}$. \\
\ \\
\noindent Let $(S, \partial_{S})$ be the complex with identical maps as in the lower picture, but with the black circle representing $u_{1}$, the white circle representing $u_{2}$, and the gray circle representing $u_{3}$. We obtain a complex over $\F(u_{1},u_{2},u_{3}, z_{7},\ldots, z_{n})$. This complex is stably chain to both $R_{1}$ and $R_{2}$. In fact $R_{1} \cong S \otimes \field{T}$ under the inclusion $u_{1} \rightarrow y_{1} + y_{4}$, $u_{2} \rightarrow y_{2} + y_{5}$, and $u_{3} \rightarrow y_{3} + y_{6}$, while $R_{2} \cong S \otimes \field{T'}$ using $u_{1} \rightarrow x_{3} + x_{4}$, $u_{2} \rightarrow x_{1} + x_{5}$, and $u_{3} \rightarrow x_{2} + x_{6}$. Thus $R_{1}$ is stably chain isomorphic to $R_{2}$.\\
\ \\
\noindent Furthermore, $n_{+}(T) = n_{+}(T')$, so 
$$
\underline{SK}(T) \simeq R_{1}[-n_{+}(T)] \cong S[-n_{+}(T)] = S[-n_{+}(T')] \cong R_{2}[-n_{+}(T')] \simeq \underline{SK}(T')
$$
Similar arguments hold for any other RIII move, so the stable chain homotopy type of $\underline{SK}(T)$ is preserved by RIII moves. 

\section{The $\delta$-graded complex for alternating tangle diagrams}\label{sec:alternating}

\noindent Let $\mathcal{T}$ be a tangle in $(\Sigma \times (-1,1), B \times \{0\})$ which admits an alternating diagram in $(\Sigma, B)$ that can be checkerboard colored. As in Lee's theorem for Khovanov homology, \cite{Lee}, the homology groups $HT_{\ast}(\mathcal{T}, b)$ can be described explicitly.\\
\ \\
\noindent First, a checkerboard coloring allows us to further decompose $CT_{\ast}(T,b)$. For any resolution $r$ the coloring identifies a subsurface of black regions $B_{r}$ in $r(T)$. Furthermore, the coloring descends to the regions of $\gly{r}$, where the color change across circles assigned an odd label in $\Z$. 

\begin{lemma} \label{lemma:2} Let $T$ be a tangle diagram for $(\Sigma, B)$ equipped with a checkerboard coloring. Let $b \in \glyph{\Sigma}{B}$ and $r \in \mathcal{O}_{n}(T,b)$. If $\pairing{\partial_{n} r}{r'} \neq 0$ then $\chi(B_{r}) = \chi(B_{r'})$. 
\end{lemma}

\noindent{\bf Proof:} From Figure \ref{fig:pic13}, $r'$ is obtained from $r$ by changing the resolution at two crossings $c_1$ and $c_2$. Changing the resolution solely at $c_1$ will either 1) split a noncontractible circle $N$ in $r(T)$ into a noncontractible circle $N'$ and a contractible circle $C$ (as in cases $II(b),II(c),IV(b)$) or 2) merge two noncontractible circles across the annulus $A$ they cobound, resulting in a contractible circle $C$ (as in cases $III(a),III(b),IV(c)$). We call the black regions in this state $B'$. If $C$ bounds a black region then $\chi(B_{r}) = \chi(B') + 1 - 2$ since we can obtain $B_{r}$ by gluing a black square to $B'$ along two edges. Thus $\chi(B_{r}) = \chi(B') - 1$. To obtain $B_{r'}$ we change the resolution at $c_{2}$. This must alter $C$ by either merging it with a noncontractible circle or by joining it to itself to get an annulus. In either case, we glue a black square to $B'$ along two edges. Thus $\chi(B_{r'}) = \chi(B') - 1$ as well. If $C$ bounds a white region, then we have glued a black square to $B_{r}$ to obtain $B'$ so $\chi(B') = \chi(B_{r}) - 1$. However, we also need to glue a white square to $B_{r'}$ to obtain $B'$, so $\chi(B_{r}) = \chi(B_{r'})$. $\Diamond$\\
\ \\
\noindent Let $\Sigma'$ be the union of the black regions in $\gly{r}$. The proof actually shows that $\chi(B_{r}) = \chi(\Sigma') = \chi(B_{r'})$.
This follows from $\gly{r'} = \gly{r}$, and  the effect that $\partial^{b}_{L}$ has on the black regions. That effect is to divide a black (or white) disc from $B_{r}$ (or $W_{r}$) and reglue it, or to divide or introduce an monocolor annulus. Since annuli are collapsed in the passage to $\gly{r}$ and have Euler characteristic $0$, $\chi(B_{r}) = \chi(\Sigma')$. 

\ \\
\ \\
\begin{defn}
A checkerboard coloring for a glyph $g \in \glyph{\Sigma}{B}$ is a choice of colors in $\{B,W\}$ for each region in the complement of the circles and arcs in $g$, satisfying the restriction that regions abutting a circle labeled with an odd integer, or an arc, have different colors, while regions abutting a circle labeled with an even integer receive the same color. 
\end{defn}

\noindent By Lemma \ref{lemma:2}, given a checkerboard coloring on $T$ we can decompose our $\delta-$graded complex along checkerboard colored glyphs $b^{c} \in \gly{r}$: 

$$
\lk O_{nc,i}(T, b^{c}) = \big\{r \big| r \in  \lk O_n(T,b), |r|=i , r\mathrm{\ colors\ }b\mathrm{\ as\ }b^{c} \big\}
$$
and:
$$
CT_{i}(T,b^{c}):= Span_{\field{T}} \big \{ r \in \lk O_{nc,i}(T, b^{c})  \big\} 
$$
equipped with the restriction of $\partial_{n}$.  \\

\begin {theorem} \label{thm:lastone} Let $T$ be a connected, alternating, checkerboard colored, tangle diagram for $(\Sigma, B)$. For each colored glyph $b^{c}$ with $b \in \glyph{\Sigma}{B}$ the chain complex $CT_{*}(T,b^{c})= \bigoplus_{i \in Z} CT_{i}(T,b^{c})$ is supported in a single homological grading when it is non-trivial. Consequently, $\underline{HT}_{\ast}(\mathcal{T}, b^{c}) \cong CT_{*}(T,b^{c})$ and any non-trivial homology groups occur in a single grading.
\end {theorem}

\noindent{\bf Proof:} Since $T$ is alternating and connected, we can switch the colorings of regions, if necessary, to ensure that the $1$-resolution at each crossing $c$ connects the two black quadrants at $c$. Let $R_{1}, \ldots, R_{k}$ be the black regions in all $0$ resolution of $T$. Then $\chi(B_{r}) = \sum_{i=1}^{k} \chi(R_{i}) - \#(\mathrm{crossings\ }c{\ with\ }r(c)=1)$ since we glue the black regions $R_{i}$ along an interval for each one resolution. Thus the homological grading $|r| = \sum_{i=1}^{k} \chi(R_{i}) - \chi(B_{r}) = \sum_{i=1}^{k} \chi(R_{i}) - \chi(B_{\gly{r}})$. The first of these terms depends only on the checkerboard coloring of $L$ while the second depends only on $b^{c}$. Thus, given the coloring on $L$ and a coloring on $b$, every $r$ representing that coloring on $b$ has resolutions in the same homological grading. The quantum grading is
$|r| + k(r)$, so the $\delta = 2|r| - (|r| + k(r)) = |r| - k(r)$ depends only on $b^{c}$. $\Diamond$.\\
\ \\
\noindent In section \ref{sec:example} there is an example of an alternating link on a genus two surface with resolutions in distinct homological gradings, but which correspond to the same uncolored glyph. Including the checkerboard coloring, however, distinguishes the two sets of generators. Furthermore, the example illustrates one of the implications of \ref{thm:lastone}: $|r| + \chi(B_{\gly{r}})$ is constant regardless of the state $r$. \\
\ \\
\noindent However, in some cases we will not need the Euler characteristic of the black region to be explicitly provided. For example, let $\Sigma$ be $\R^{2}$ minus $n$ open discs indexed in some way by $\{1,\ldots, n\}$ and let $B = \emptyset$. We will consider a connected alternating link diagram $L$ in $\Sigma$. Taking the mirror of $L$ if necessary, we can arrange both that $1$-resolutions join black regions, and that the unbounded region of $\Sigma \backslash T$ is colored black. There is then a unique checkerboard coloring for any glyph $g \in \glyph{\Sigma}{\emptyset}$. Consequently, the
homological grading depends only on $g$.   \\
\ \\
\noindent As an example of the latter case, consider $L$ in $\R^{2} \backslash\{(0,0)\}$. There is a glyph $g_{n}$ for each $n \in \Z$, represented by the annulus  with one circle labeled by $n$. The Euler characteristic of the black regions will be $0$, and  the $\delta$-grading will be determined
by the glyph $g_{n}$ alone, with $k = n$. More generally, 

\begin{prop}
Let $\Sigma_{n+1} = \R^{2} \backslash\{(0,0),(1,0), \ldots, (n,0)\}$. Let $\mathcal{L} \subset \Sigma \times \R$ be a link with a connected, alternating projection $L \subset \Sigma$, checkerboard colored as above. Let $m_{B}$ be the number of black regions surrounding a puncture. Let $J = 0$ if the unbounded region is colored black, and $1$ if the unbounded region is colored white. For any $g \in \glyph{\Sigma}{\emptyset}$, if $\underline{HT}_{\ast}(L,g) \neq 0$ then it is supported in the grading $$\delta=  \sigma(\mathcal{L}) - k(g) - \chi(B_{g}) - m_{B} + J$$ where $\sigma(\mathcal{L})$ is the signature of the link in $\R^{3}$ found by applying $\Sigma \times \R \hookrightarrow \R^{2} \times \R$ to $\mathcal{L}$ . 
\end{prop} 

\noindent {\bf Proof:} For a diagram consisting of noncontractible circles the quantum grading of $r$ equals $|r| + k(r)$. Thus $\widetilde{\delta} = 2|r| - (|r| + k(r)) = |r| - k(r)$. If $r$ represents $g$ then $k(r) = k(g)$, so $\widetilde{\delta}(r) = |r| - k(g)$. Let $\#(B)$ be the number of black regions. Since the diagram is connected, the black regions, ignoring punctures, consist of discs if the unbounded region is colored white, or discs and a single annulus if the unbounded region is colored black. To compute $\chi(B_{g})$ we glue the black regions at $1$-resolved crossings, and puncture the $m_{B}$ regions surrounding an $(i,0)$ point. Thus $\chi(B_{g}) = \#(B) - I - |r| - m_{B}$ where 1) $I=0$ if the unbounded region is white, and $1$ otherwise, 2) we subtract $1$ for each gluing at a $1$-resolution, and 3) we subtract $1$ for each puncture in a black region. Thus $\widetilde{\delta}(r) = \#(B) - I - \chi(B_{g}) - m_{B} - k(g)$. To obtain the corresponding generator in $\underline{HT}_{\ast}(L,g)$ we must subtract $n_{+}(L)$ to obtain $\delta = \#(B) - n_{+} - I - \chi(B_{g}) - m_{B} - k(g)$. However, $\sigma(\mathcal{L}) = \#(B) - 1 - n_{+}(L)$ by the formula of Gordon-Litherland. Thus, $\delta = \sigma(\mathcal{L}) + 1 - I - \chi(B_{g}) - m_{B} - k(g)$. $\Diamond$\\
\ \\
\noindent For the annulus $\Sigma_{1}$, the formula above simplifies, and confirms the result for the Khovanov skein homology of an alternating, annular link proved by the second author in \cite{LR2}. 

\begin{cor}
Let $L$ be a connected, alternating link diagram in $\R^{2}\backslash\{(0,0)\}$ checkerboard colored so that the $1$-resolutions merge the black regions, and let $\mathcal{L}$ be the corresponding link in $\R^{3}$. Let $g_{k}$ be the glyph determined by the unit circle labeled with $k \in \Z$. Then any non-trivial groups $\underline{HT}(L,g_{k})$ are supported in the $\delta$-gradings determined by
\begin{enumerate}
\item if $k$ is odd, then $\delta + k = \sigma(\mathcal{L})$,
\item if $k$ is even and the unbounded region is colored black, then $\delta + k = \sigma(\mathcal{L}) - 1$,
\item while if $k$ is even and the unbounded region is colored white, then $\delta + k = \sigma(\mathcal{L}) + 1$
\end{enumerate}
\end{cor} 

\noindent{\bf Proof:} If $k$ is odd then either $C$'s interior or its exterior is colored black, but not both. In the first case $J=1$, $m_{B}=1$ and $\chi(B_{g})=0$, so $\delta + k = \sigma(\mathcal{L})$. The same formula holds in the second case: $J=0$, $m_{B}=0$ and $\chi(B_{g})=0$ so $\delta + k = \sigma(\mathcal{L})$. On the other hand, if $k$ is even then both interior and exterior will be colored the same color. If they are colored black then $J=0$, $m_{B}=1$ and $\chi(B_{g})=0$ so $\delta + k = \sigma(\mathcal{L}) - 1$. If they are both colored white then $J=1$, $m_{B} = 0$ and $\chi(B_{g})= 0$, so $\delta + k = \sigma(\mathcal{L}) + 1$. $\Diamond$\\
\ \\
\noindent The next simplest case is that of tangles in $D^{2}$ with ends on $2n$ points in the boundary. We consider $D^{2}$ embedded in $S^{2}$, and use a crossingless matching of the $2n$-points in $D^{2}$ to complete the tangle diagram to a link diagram. 

\begin{cor}
Let $T$ be a connected, oriented, alternating, checkerboard colored tangle diagram in $(D^{2}, B_{2n})$ whose regions are homeomorphic to discs. Let $g$ be a colored glyph in $\glyph{D^2}{B_{2n}}$ which colors the same boundary arcs black as does $T$. Let $b$ be any crossingless matching of $B_{2n}$ which can be glued to $T$ to give the diagram of an {\em oriented} link $\mathcal{L}$. Let $g \# b$ be the colored decomposition of $S^{2}$ found by gluing the arcs of $g$ to those of $b$. Then $\underline{HT}(T,g)$ is either trivial or supported in the grading 
$$
\delta = \sigma(\mathcal{L}) - \chi(B_{g\# b}) + 1
$$
where $\#(B_{g \#g_{T}})$ is the number of black regions in $g \# b$.
\end{cor}

\noindent{\bf Note:} The glyph $g$ is just a crossingless matching, as any circle in $D^{2}$ is contractible. Also, in this case $\underline{HT}(T,g)$ will be homeomorphic to the untwisted skein homology, so this result also establishes a Lee type theorem in the context of \cite{APS2}.\\
\ \\
\noindent{\bf Proof:} Each black region in $T$ is a disc, so a resolution will have glyph $g$ if and only if the $1$-resolutions are chosen to merge the discs into black regions homeomorphic to those in $g$. For each merge we will join two discs, and thus the number of $1$ resolutions will be $\#(B_{T}) - \#(B_{g})$ where $B_{T}$ is the number of black regions in $T$. After shifting by $-n_{+}(T)$ we have generators in $\delta - k(g) =  \#(B_{T}) - \#(B_{g}) - n_{+}(T)$. Since there are no noncontractible circles, $k(g) = 0$ and $\delta = \#(B_{T}) - \#(B_{g}) - n_{+}(T)$. $\sigma(\mathcal{L}) = \#(B_{L}) - n_{+}(T) - 1$, so $\delta = \#(B_{T}) - \#(B_{g}) + \sigma(\mathcal{L}_{T}) - \#(B_{L}) + 1 $ $= \sigma(\mathcal{L}_{T}) + 1 - \#(B_{g}) + \#(B_{T}) - \#(B_{L})$. However, each black region in the interior of $D^{2}$ contributes to both $B_{T}$ and $B_{L}$. There are $n$ more discs in $B_{T}$, so $\#(B_{T}) - \#(B_{L}) = n - \#(\mathrm{\ black\ regions\ in\ }L\mathrm{\ intersecting\ }S^{1})$. The latter is just $\#(B_{b})$ since we just join discs to the edges of the regions in $B_{b}$ to obtain the black regions in $L$ intersecting $S^{1}$. $B_{b}$ and $B_{g}$ consist of unions of discs, and their union along the $n$-arcs in the boundary has Euler characterisitic $\chi(B_{g\#b}) = \#(B_{b}) + \#(B_{g}) - n$. This gives the formula above $\Diamond$.\\
\ \\
\noindent It remains to see that there is always a link $\mathcal{L}$ as in the corollary. We give a construction of the requisite matching in the proof of

\begin{lemma}
Let $T$ be a connected, alternating, oriented tangle diagram in $(D^{2}, B)$ where $B$ is the set of $2n^{th}$ roots of unity in $S^{1}$. If the regions of $T$ in $D^{2}$ are all homeomorphic to discs, there is an alternating, oriented link diagram $L_{T}$ in $S^{2}$ whose intersection with $D^{2}$ is precisely $T$, and which has no additional crossings. 
\end{lemma}

\noindent{\bf Proof:} We label each end of $T$ on $S^{1}$ with a pair from $\{+,-\}\times\{u,o\}$ where a $+$ indicates that the tangle points out of $D^{2}$ at that end and a $-$ means that the tangle points in. A $u$ indicates that the strand should go under at the {\em next} crossing when heading out of the disc, if one wishes to maintain alternation. We need to find arcs embedded in the complement of $D^{2}$ which join $+$ ends to $-$ ends and $u$ ends to $o$ ends. In other words, if one end of the arc is labeled $(+,u)$ then the other end needs to be labeled $(-,o)$. We know that there are $n$ ends labeled with each of $+$ and $-$, but that there is otherwise no pattern. The pattern of $u$'s and $o$'s, on the other hand, can be determined: going around $S^{1}$ counter-clockwise will yield an alternating pattern of $u$'s and $o$'s. To see this, take the arc, $a$, going from an end $e_{i}$ marked with a $u$ to the end $e_{i+1}$. By assumption, this is in the boundary of a region $R$ which is homeomorphic to a disc. The arc $a$ and the arcs in $T$ which abut $R$ make $R$ into a polygon. Going from $e_{i}$ to $e_{i+1}$ around the path in the polygon we come first to an overstrand, and then turn onto the understrand that is the next arc. We repeat this some number of times and then turn onto the understrand which leads to $e_{i+1}$. Since we came to $e_{i+1}$ from an understrand, to have alternation requires an overstrand at the next crossing when heading out of the disc. Thus $e_{i+1}$ will be an $o$-end. A similar argument applies to $e_{i}$ as an $o$-end to see that $e_{i+1}$ is a $u$-end. To find $L_{T}$ start join each $+$ end at $e_{i}$ to $e_{i+1}$ if it is a $-$ end. Push this arc into the disc so that we have fewer ends, and then repeat. The result is an alternating tangle, so the labeling will preserve the $o-u$ alternation. After finitely many repetitions we will have joined all the ends together. Due to the coloring convention we will color black the regions abutting arcs from $e_{i}$ with an $o$ to $e_{i+1}$ with a $u$. $\Diamond$


\begin{thebibliography}{99}

\bibitem{APS1} M. Asaeda, J. Przytycki, A. Sikora, {\em Categorification of the Kauffman bracket skein module of I-bundles over surfaces}. Algebr. Geom. Topol. 4 (2004), 1177–1210 (electronic).

\bibitem{APS2} M. Asaeda, J. Przytycki, A. Sikora, {\em Categorification of the skein module of tangles}. Primes and knots, 1–8,
Contemp. Math., 416, Amer. Math. Soc., Providence, RI, 2006.

\bibitem{Bald}J. Baldwin, A. Levine, {\em A combinatorial spanning tree model for knot Floer homology}. arXiv: math.GT/1105.5199

\bibitem{Bar1}D. Bar-Natan, {\em On Khovanov's categorification of the Jones polynomial}. {\em Alg. \& Geom. Top.} 2:337--370 (2002).

\bibitem{Bar2}D. Bar-Natan, {\em Khovanov's homology for tangles and cobordisms}. {\em Geom. Topol.} 9:1443-1499 (2005).


\bibitem{TJ}T. Jaeger, {\em A Remark on Roberts' Totally Twisted Khovanov Homology}.  arXiv:math.GT/1109.1805 

\bibitem{Khov}M. Khovanov, {\em A categorification of the Jones polynomial}. {\em Duke Math. J.}  101(3):359--426 (2000). 

\bibitem{Kriz}D. Kriz, I. Kriz {\em Baldwin-Ozsváth-Szabó cohomology is a link invariant}. arXiv:math.GT/1109.0064

\bibitem{Lee}E. S. Lee, {\em An endomorphism of the Khovanov invariant}. Adv. Math. 197(2):554-–586 (2005).

\bibitem{Mann}A. Manion, {\em A sign assignment in totally twisted Khovanov homology}. arXiv:math.GT/1109.4982


\bibitem{LR}L. Roberts, {\em Totally Twisted Khovanov Homology}. arXiv:math.GT/1109.0508.

\bibitem{LR2}L. Roberts, {\em On knot Floer homology in double branched covers}. arXiv:math.GT/0706.0741.

\bibitem{Viro}O. Viro, {\em Khovanov homology, its definition and ramifications}. {\em Fund. Math.} 184:317--342 (2004).


\end{thebibliography}
\end{document}